\documentclass[letterpaper]{mc2025}

\usepackage{color}
\usepackage{subcaption}
\usepackage{bm}                 
\usepackage{float}
\usepackage{comment}

\newcommand{\dx}{\Delta x}
\newcommand{\dxi}{\Delta \xi}
\newcommand\eref[1]{Eq.~(\ref{#1})}
\newcommand\fref[1]{Figure~\ref{#1}}
\newcommand\tref[1]{Table~\ref{#1}}

\usepackage{nomencl} 
\makenomenclature

\title{CWENO Interpolation for Non-Oscillatory Stochastic Collocation in~Uncertainty Quantification Problems}

\addAuthor{Alina Chertock}{1}
\addAuthor[\href{mailto:aiskhak@ksu.edu}{aiskhak@ksu.edu}]{Arsen S. Iskhakov}{2}
\addAuthor{Anna Iskhakova}{2}
\addAuthor{Alexander Kurganov}{3}

\addAffiliation{1}{Department of Mathematics, North Carolina State University, Raleigh, NC, USA}
\addAffiliation{2}{Department of Mechanical and Nuclear Engineering, Kansas State University, Manhattan, KS, USA}
\addAffiliation{3}{\mbox{Department of Mathematics, Shenzhen International Center for Mathematics, and Guangdong}\newline\mbox{Provincial Key Laboratory of Computational Science and Material Design, Southern University}\newline\mbox{of Science and Technology, Shenzhen 518055, Guangdong, China}}

\Abstract{Uncertainty quantification (UQ) in mathematical models is essential for accurately predicting system behavior under variability. This study provides guidance on method selection for reliable UQ across varied functional behaviors in engineering applications. Specifically, we compare several interpolation and approximation methods within a stochastic collocation (SC) framework, namely: generalized polynomial chaos (gPC), B-splines, shape-preserving (SP) splines, and central weighted essentially nonoscillatory (CWENO) interpolation, to reconstruct probability density functions (PDFs) and estimate statistical moments. These methods are assessed for both smooth and discontinuous functions, as well as for the solution of the 1-D Euler and shallow water equations. While gPC and interpolation B-splines perform well with smooth data, they produce oscillations near discontinuities. Approximation B-splines and SP splines, while avoiding oscillations, converge more slowly. In contrast, CWENO interpolation demonstrates high robustness, effectively capturing sharp gradients without oscillations, making it suitable for complex, discontinuous data. Overall, CWENO interpolation emerges as a versatile and effective approach for SC, particularly in handling discontinuities in UQ.}


\shortTitle{M&C 2025 Full Summary Template}

\authorHead{A. Chertock, A.S. Iskhakov, A. Iskhakova, A. Kurganov}

\begin{document}
\section{Introduction}
Many engineering problems involve inherent uncertainties arising from various sources. Uncertainty quantification (UQ), whether in input parameters or in initial and boundary conditions (BCs), often due to empirical approximations or measurement errors, is crucial to conducting sensitivity analysis and gaining insight to improve model accuracy. This study explores different stochastic collocation (SC) methods for UQ to address such challenges. Consider a mathematical model where such uncertainties can be represented using random variables:
\begin{equation}
\bm{U}_t=\mathcal{L}(x,t,\bm{U},\bm{U}_x,\ldots;\xi)
\label{eq:1.1.1}
\end{equation}
where $x$ is the spatial variable, $t$ is the time, $\bm{U}(x,t;\xi)$ is an unknown vector-function, $\mathcal{L}$ is a differential operator, $\xi$ is a real-valued random variable with known probability density function (PDF) $p_\xi(\xi)$.

Among the various UQ methods, Monte Carlo approaches~\cite{mishra2013} are known for their reliability but are computationally expensive due to the large number of realizations needed. A more efficient alternative is SC, where the governing equations are first solved using a deterministic numerical solver at a set of collocation points in random space. Using the solution obtained, a surrogate model is constructed using interpolation or approximation techniques. This model enables efficient prediction of the solution at additional sample points with subsequent reconstruction of the PDF and estimation of statistical moments~\cite{xiu2009}. In building the surrogate, several interpolation or approximation techniques are available. In this work, we compare the performance of generalized polynomial chaos (gPC), B-splines, shape-preserving (SP) splines, and CWENO (central weighted essentially non-oscillatory) interpolation methods.

In gPC, the solution is represented as a series of orthogonal polynomials~\cite{maitre2010}. The application of gPC methods to~\eref{eq:1.1.1} presents specific challenges, particularly in handling shock waves, which frequently arise when solving nonlinear hyperbolic PDEs. Although these discontinuities occur in the spatial domain, their propagation speed is affected by uncertainty, leading to discontinuities in the random variable and resulting in Gibbs-type phenomena~\cite{maitre2004}. In our previous work~\cite{chertock2024}, to avoid oscillations, we evaluated the performance of B-splines~\cite{deboor1972} and SP rational quartic splines~\cite{zhu2015} as approximation/interpolation tools in SC. Although approximation B-splines avoid oscillations near discontinuities, their tendency to smear data limits their effectiveness for accurate UQ. In contrast, SP splines apply only localized smoothing, preserving solution accuracy in other regions and maintaining the asymptotic behavior of the solution~\cite{dai2021}.

In this paper, we extend our investigation to compare the capabilities of gPC, B-splines, SP splines, and CWENO interpolation for UQ by reconstructing the PDF and estimating statistical moments. Although both B-splines and SP splines effectively minimize oscillations, they have limited accuracy and convergence because of inherent smoothing effects. This limitation highlights the need for alternative approaches. We propose the use of CWENO interpolation methods~\cite{semplice2021}, which achieve high accuracy in smooth regions while avoiding oscillations near discontinuities. Thus, CWENO interpolation offers an effective and adaptable approach for SC, particularly suited for complex, discontinuous data.

\section{Methodology}
We begin by selecting a set of collocation points in the random space, denoted as $\{\xi_n\}_{n=1}^{N}$, and numerically solving the following deterministic systems for each point:
\begin{equation*}
(\bm{U}(x,t;\xi_n))_t=\mathcal{L}(x,t,\bm{U}(x,t;\xi_n),(\bm{U}(x,t;\xi_n))_x,\ldots;\xi_n), \quad n=1,\ldots,N,
\end{equation*}
up to a specified final time $T$. For each discrete spatial node, we perform an interpolation/approximation in the random space using gPC, B-splines, SP splines, and CWENO. The interpolation/approximation is carried out at points $\{\xi_m\}_{m=1}^{M}$, $M \gg N$, which are randomly sampled according to the known $p_\xi(\xi)$. For each component $U$ of $\bm{U}$, the resulting values, $U_m=U(x,T;\xi_m)$, are then used to reconstruct the PDF, $p(U)$, via the histogram method and to estimate the mean and variance:
\begin{equation*}
\mathbb{E}[U]=\sum_{i=1}^{n_\text{bins}} p_i U_i \Delta U, \quad \mathbb{V}[U]=\sum_{i=1}^{n_\text{bins}} p_i (U_i - \mathbb{E}[U])^2 \Delta U,
\end{equation*}
where $p_i=p(U_i)$ is the reconstructed PDF value  within the $i^\text{th}$ histogram bin, $U_i$ is the midpoint value within the $i^\text{th}$ bin, and $\Delta U$ is the width of the histogram bins, calculated as $\Delta U=\frac{U_\text{max}-U_\text{min}}{n_\text{bins}+1}$. The number of bins, $n_\text{bins}$, is determined by the \texttt{`auto'} option in the \texttt{numpy.histogram} function available in \texttt{Python}.

While we omit detailed descriptions of the gPC, B-splines, SP splines, and CWENO interpolation methods (can be found in the cited literature), it is important to highlight the following key points:
\begin{itemize}
    \item In gPC~\cite{maitre2010}, the quadrature rule is determined by $p_\xi(\xi)$. In this work, we use the Gauss-Legendre quadrature when $p_\xi(\xi)$ is uniformly distributed and the Gauss-Hermite quadrature when $p_\xi(\xi)$ follows a normal distribution. Additionally, in gPC, the mean and standard deviation estimates are readily obtained from the expansion coefficients.
    \item In B-splines~\cite{deboor1972}, the selection of knots can lead to approximation (e.g., when knots are uniformly spaced) or interpolation (when knots are optimally selected). As demonstrated in Section~\ref{sec:3}, approximation may overly smooth the solution, while interpolation produces oscillations near discontinuities.
    \item For SP splines, the specific choice of local tension shape parameters ensures the preservation of the monotonicity, positivity, and convexity of the spline, as discussed in~\cite{zhu2015}.
    \item Various CWENO interpolation methods exist~\cite{semplice2021}. In this work, we implement the CWENO-Z of order 7 blending four $3^\text{rd}$-degree polynomials and one $6^\text{th}$-degree polynomial.
\end{itemize}

\section{Numerical Examples}\label{sec:3}
\subsection{Example 1: A Smooth Function}
We begin with a smooth function assumed to be an analytical solution of~\eref{eq:1.1.1} at time $T$:
\begin{equation}
U(\xi)=3\cos(\pi\xi), 
\label{eq:3.1.1}
\end{equation}
Our objective is to approximate the PDF of $U(\xi)$ denoted by $p(U)$, as well as provide the estimates for mean and standard deviation. We first employ the Monte Carlo method to obtain the reference PDF $\tilde{p}(U)$. For that, random samples $\{\xi_m\}_{m=1}^{M}$, $M=3\times10^7$ are generated. The corresponding values of $U_m$ are then calculated using~\eref{eq:3.1.1} and utilized in the histogram method. To implement SC, we use collocation points $\{\xi_n\}_{n=1}^{N}$, with $N=7,8,\ldots,19,20,40,60$. In gPC, these points follow a specified quadrature rule, while for other methods they are uniformly spaced. This setup enables interpolation or approximation of $U_m$ at the sample points $\{\xi_m\}_{m=1}^{M}$, followed by PDF reconstruction. For consistency in error analysis, the same number of bins, $n_\text{bins}$, is used across all methods.

\begin{figure}[ht!]
\centering 
\begin{subfigure}[t]{0.30\textwidth} 
\centering
\includegraphics[width=\textwidth]{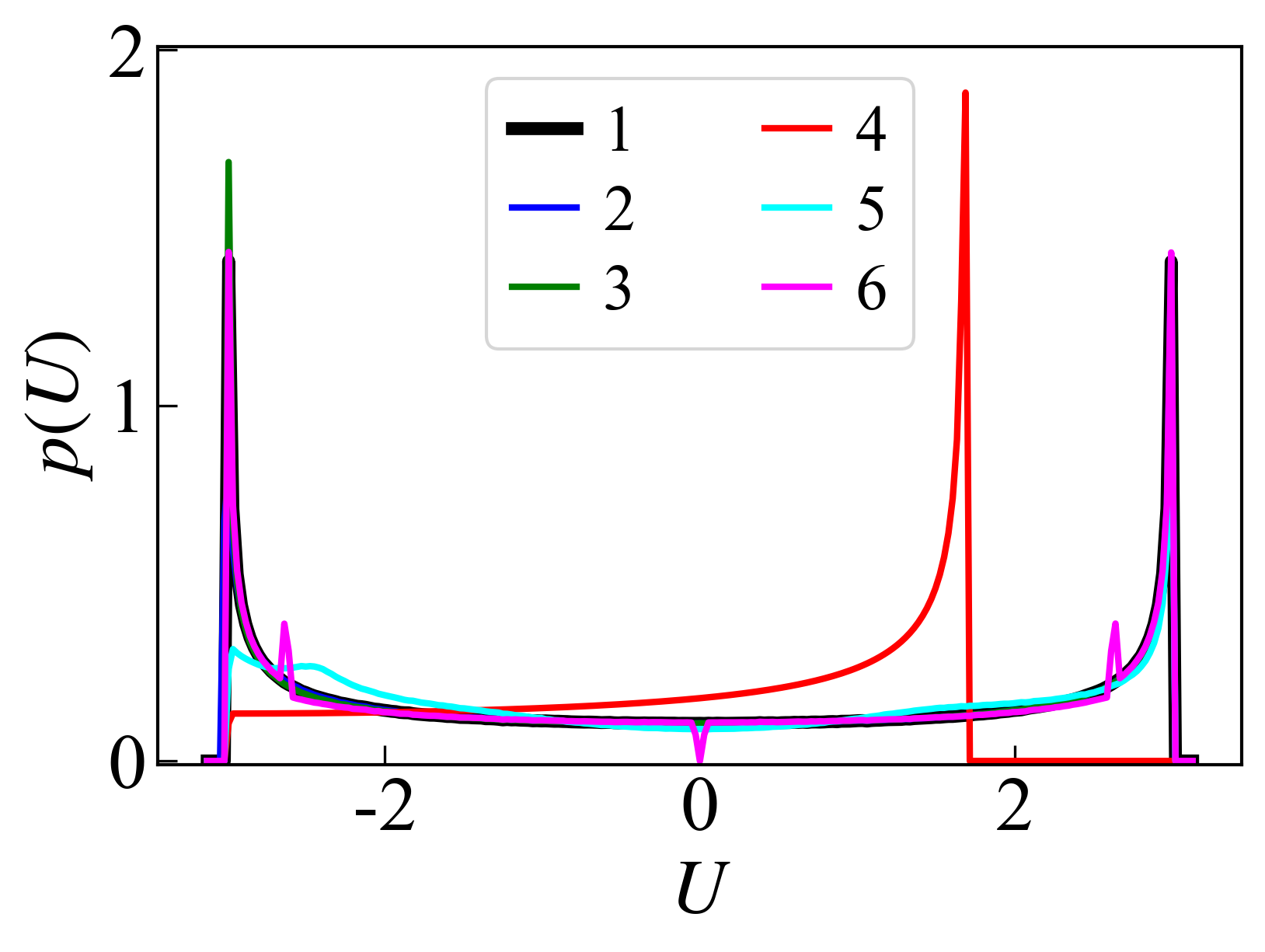}
\caption{}
\label{fig:3.1.1a}
\end{subfigure}
\begin{subfigure}[t]{0.30\textwidth}
\centering 
\includegraphics[width=\textwidth]{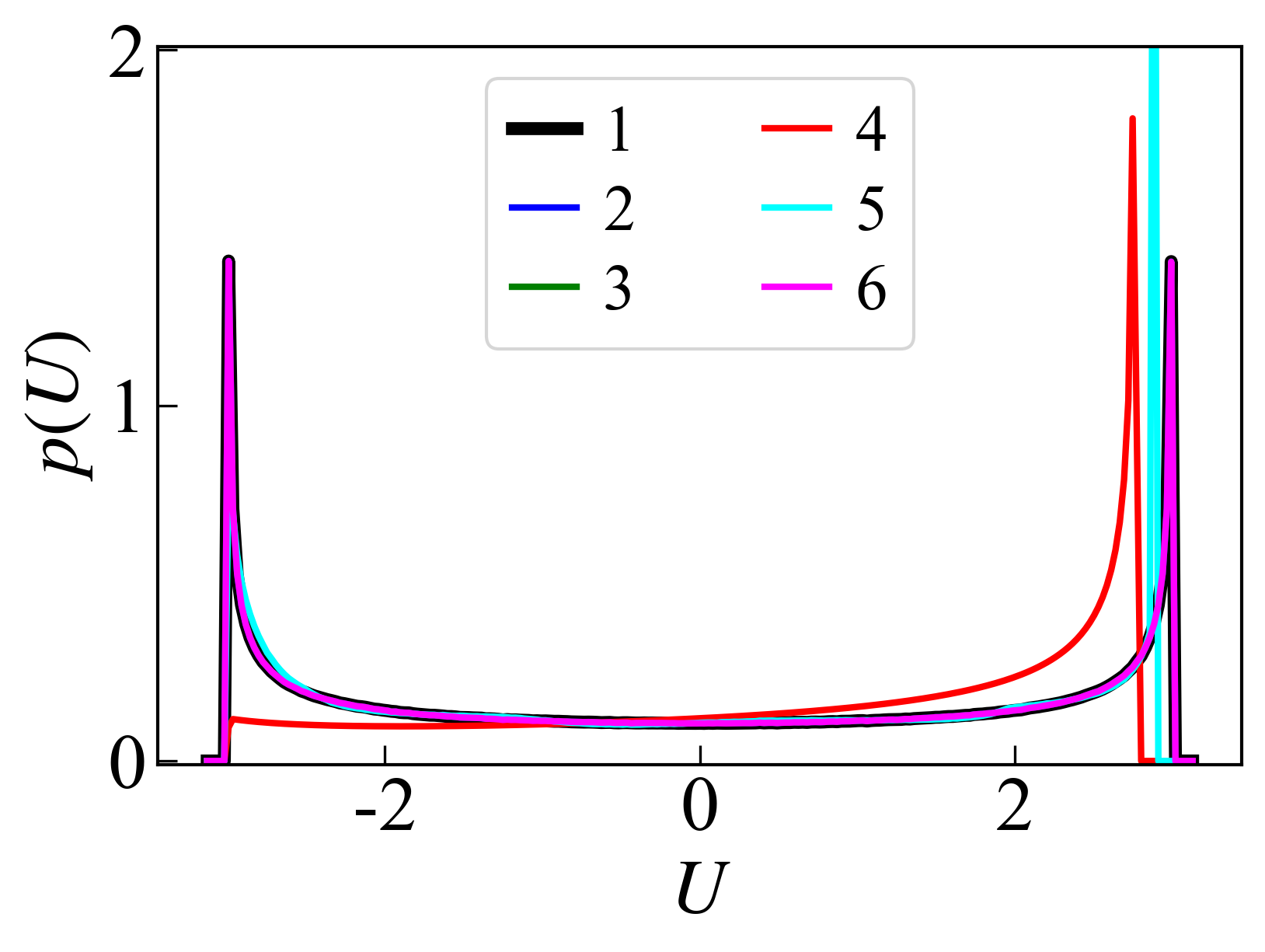}
\caption{}
\label{fig:3.1.1b}
\end{subfigure}
\begin{subfigure}[t]{0.30\textwidth}
\centering
\includegraphics[width=\textwidth]{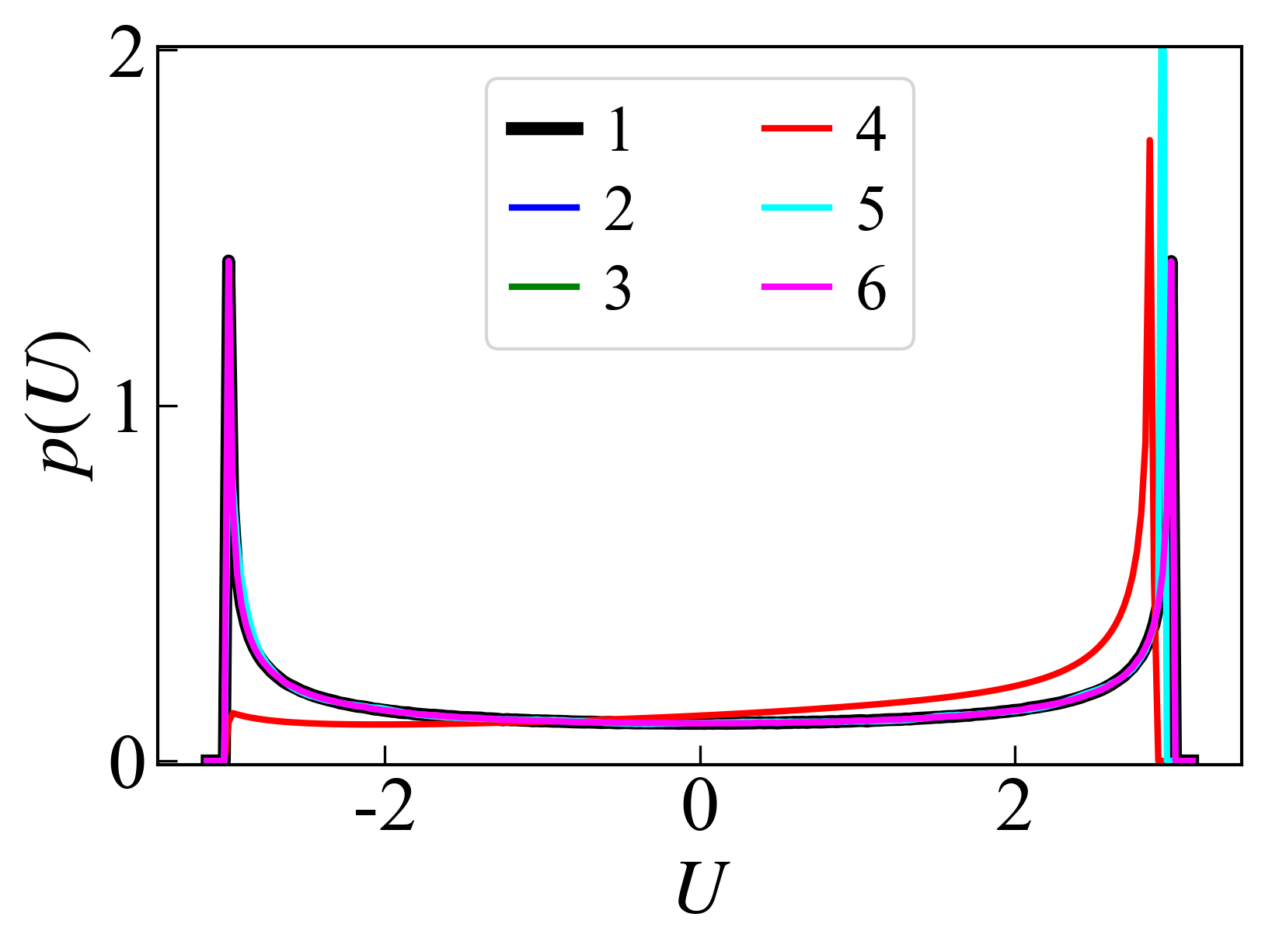} 
\caption{} 
\label{fig:3.1.1c} 
\end{subfigure} 
\captionsetup{justification=centerlast} 
\caption{Example 1, Test 1: PDF of $U(\xi)$ for (\subref{fig:3.1.1a}) $N=7$, (\subref{fig:3.1.1b}) $N=12$ and (\subref{fig:3.1.1c}) $N=16$: 1~--~reference, 2~--~gPC, 3~--~interpolation B-spline, 4~--~approximation B-spline, 5~--~SP spline, 6~--~CWENO interpolation.}
\label{fig:3.1.1} 
\end{figure}

\textbf{Test 1 -- Uniform Distribution}. In this scenario, we assume $\xi\sim\mathcal{U}[-1,1]$. \fref{fig:3.1.1} shows the PDFs obtained for $N=7,12$ and $16$. We adopt the following notation and color scheme throughout the paper: 1~(black)~--~reference, 2~(blue)~--~gPC, 3~(green)~--~interpolation B-spline, 4~(red)~--~approximation B-spline, 5~(cyan)~--~SP spline, and 6~(magenta)~--~CWENO interpolation. One can see that gPC and the interpolation B-spline provide consistently accurate results for all values of $N$. In contrast, approximation B-spline and SP spline fail to reproduce the PDF due to their smoothing properties. CWENO interpolation requires a sufficient number of collocation points $N$ to avoid abrupt jumps between the histogram intervals (can be seen for $N=7$, \fref{fig:3.1.1a}) caused by the use of local stencils.

\begin{figure}[ht!] 
\centering
\begin{minipage}[t]{0.64\textwidth}
\vspace{0pt}
\centering
\includegraphics[width=0.5\textwidth]{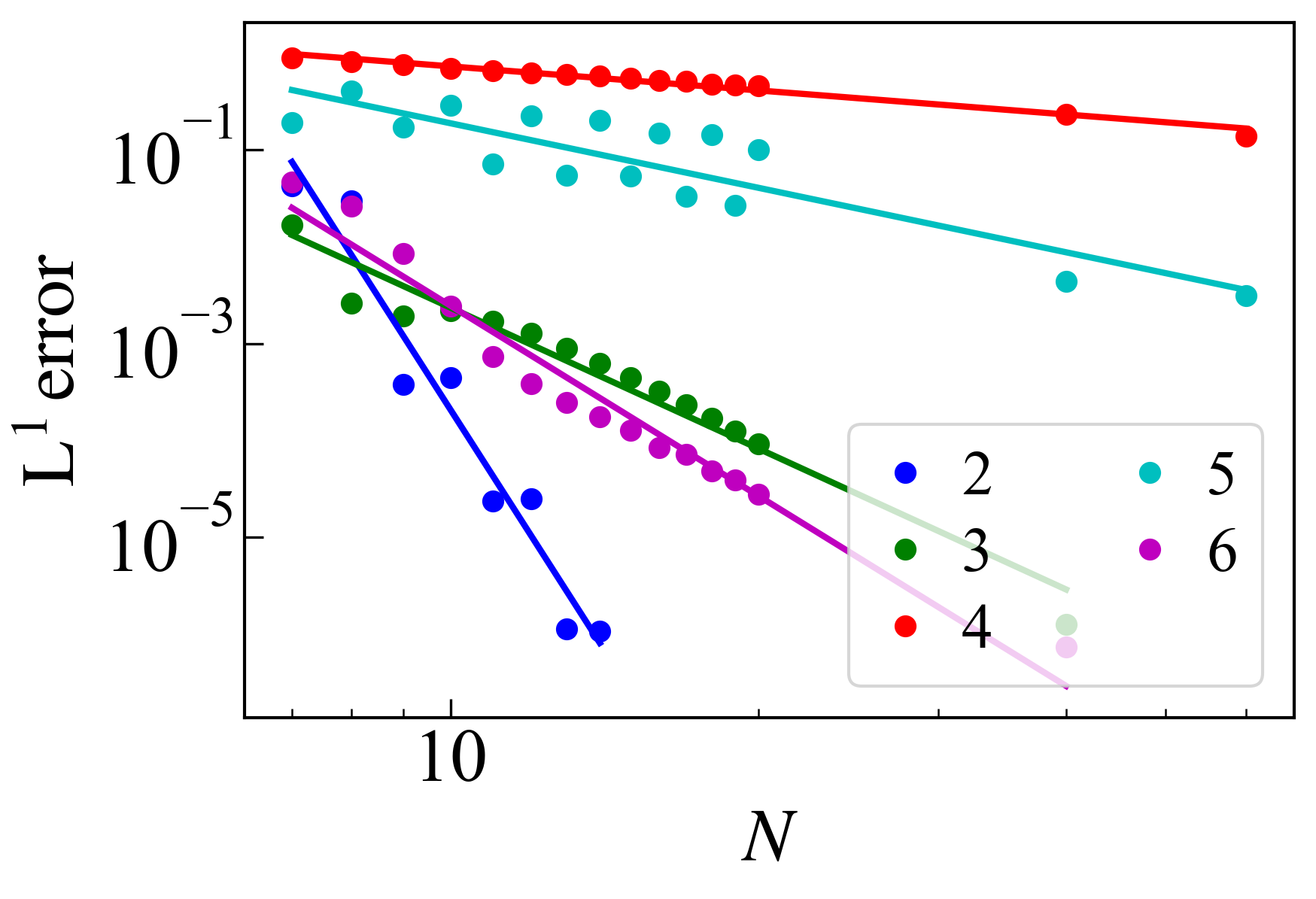}
\captionsetup{justification=centerlast} 
\caption{Example 1, Test 1: L$^1$ error of PDF reconstruction as a function of $N$, and the corresponding power-law fit (solid lines).} 
\label{fig:3.1.2}
\end{minipage}
\hfill
\begin{minipage}[t]{0.34\textwidth}
\vspace{0pt}
\centering
\captionsetup{justification=centerlast}
\captionof{table}{Example 1, Test 1: \mbox{Observed} values of $k$.}
\begin{tabular}{|l|c|}\hline
\multicolumn{1}{|c|}{\textbf{Method}} & $k$ \\\hline
gPC                     & 16.2 \\\hline
Interpolation B-spline  & 5.2  \\\hline
Approximation B-spline  & 0.8  \\\hline
SP spline               & 2.2  \\\hline
CWENO interpolation     & 6.2  \\\hline
\end{tabular}
\label{tab:3.1.1}
\end{minipage}
\end{figure}

\begin{figure}[ht!]
\centering 
\begin{subfigure}[t]{0.30\textwidth} 
\centering
\includegraphics[width=\textwidth]{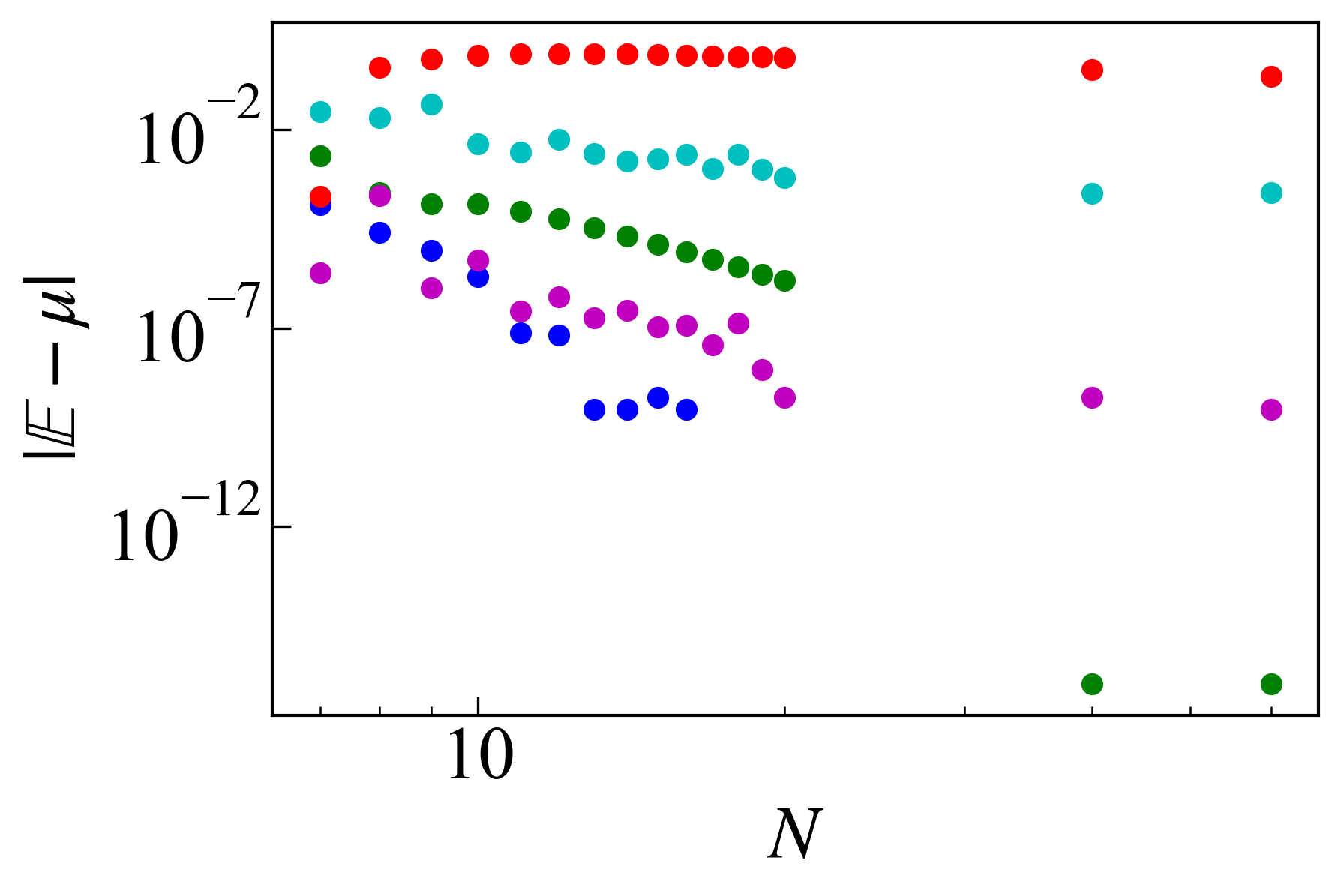}
\caption{}
\label{fig:3.1.3a}
\end{subfigure}
\begin{subfigure}[t]{0.30\textwidth}
\centering 
\includegraphics[width=\textwidth]{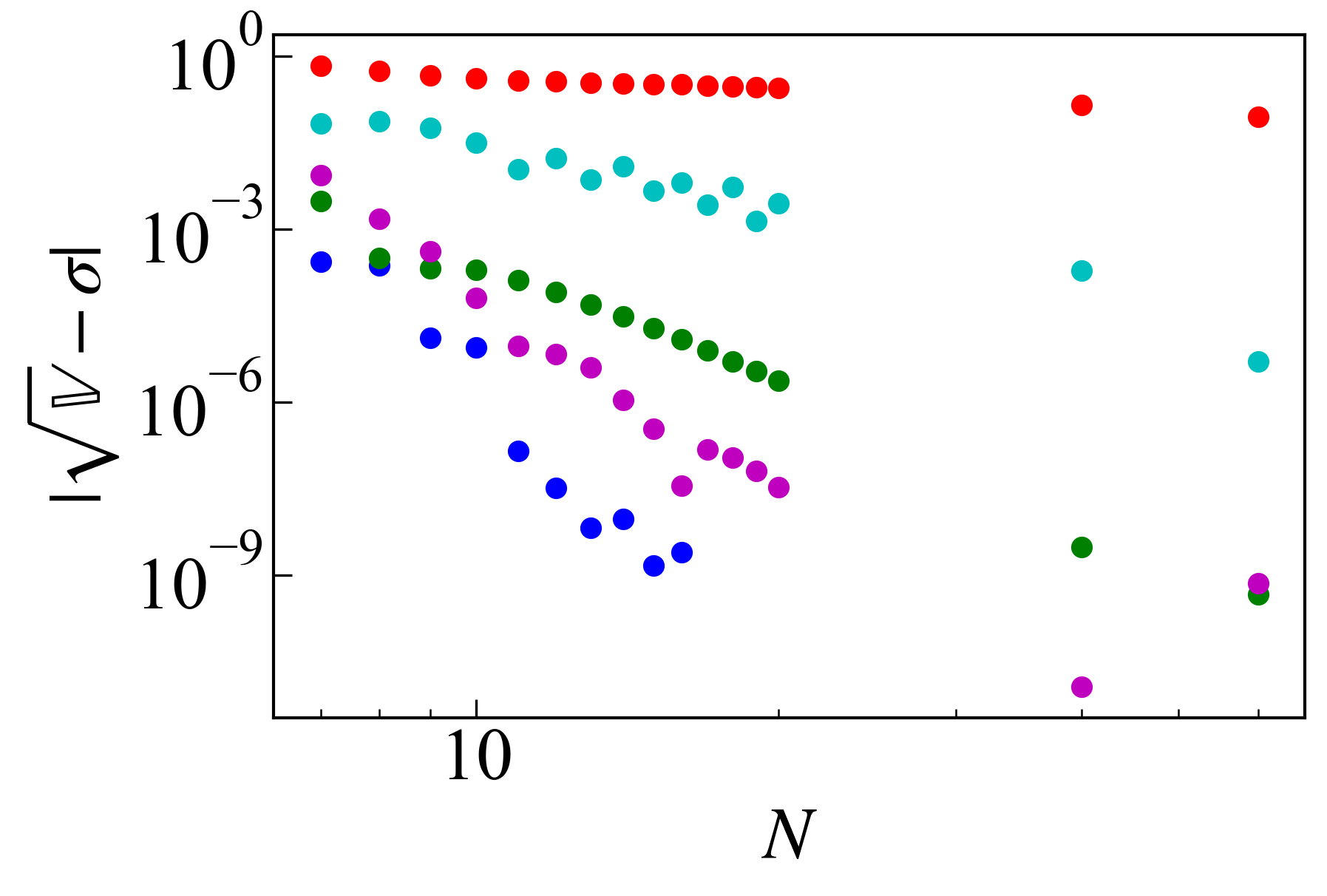}
\caption{}
\label{fig:3.1.3b}
\end{subfigure}
\captionsetup{justification=centerlast} 
\caption{Example 1, Test 1: Error of (\subref{fig:3.1.3a}) mean and (\subref{fig:3.1.3b}) standard deviation estimates for different $N$.}
\label{fig:3.1.3} 
\end{figure}

Following \cite{ditkowski2020}, we perform the convergence analysis. \fref{fig:3.1.2} shows the L$^1$ error, $||\tilde{p}(U)-p(U)||_1$, as a function of $N$. As expected, gPC exhibits the fastest convergence, followed by interpolation using CWENO and B-splines. In contrast, convergence of the approximation B-spline and SP spline is much slower, which limits their ability to produce accurate PDF in engineering problems. To further analyze the convergence rate, we apply a linear least-squares approximation:
\begin{equation*}
||\tilde{p}(U)-p(U)||_1\approx KN^{-k}
\end{equation*}
The observed values of $k$ are summarized in \tref{tab:3.1.1}. We can also analyze the convergence of the estimates for mean $\mu$ and standard deviation $\sigma$. The absolute error of the estimates $\mathbb{E}[U]$ and $\sqrt{\mathbb{V}[U]}$ are shown in~\fref{fig:3.1.3} as functions of $N$. One can see that the obtained results are, generally, in agreement with the previous discussion -- gPC, interpolation B-spline and CWENO interpolation yield more accurate predictions than the approximation B-spline and SP spline.

\textbf{Test 2 -- Normal Distribution}. Here we assume $\xi\sim\mathcal{N}(0,\sigma_\xi^2)$, with $\sigma_\xi=0.33$. This results in the Gauss-Hermite quadratures for the gPC. For the other methods, the data points are uniformly spaced within $[-6\,\sigma_\xi,6\,\sigma_\xi]$ (sample points appearing outside of this interval are re-sampled to avoid extrapolation).

\begin{figure}[ht!] 
\centering
\begin{minipage}[t]{0.64\textwidth}
\vspace{0pt}
\centering
\includegraphics[width=0.5\textwidth]{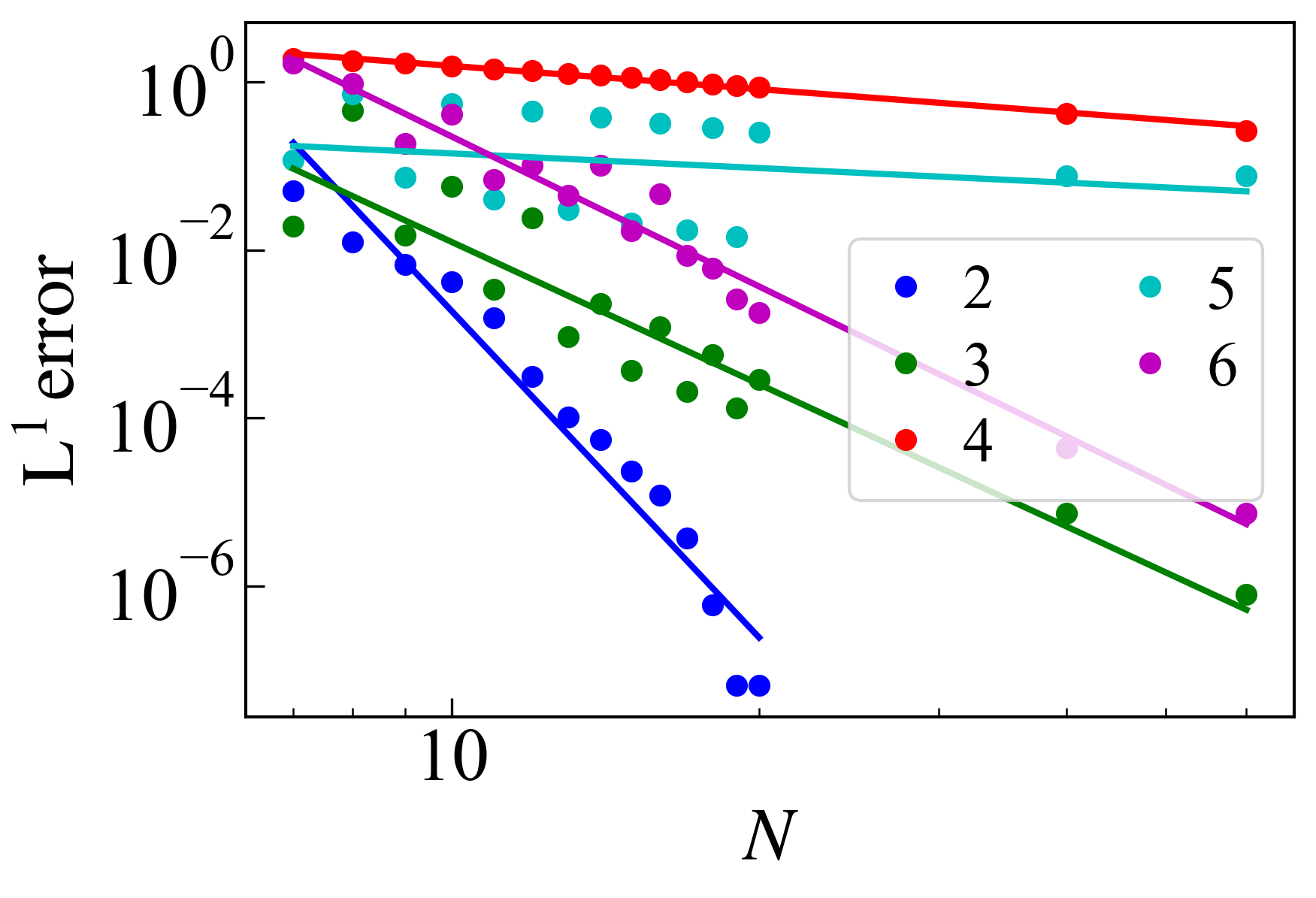}
\captionsetup{justification=centerlast} 
\caption{Example 1, Test 2: L$^1$ error of PDF reconstruction as a function of $N$, and the corresponding power-law fit (solid lines).}
\label{fig:3.1.4}
\end{minipage}
\hfill
\begin{minipage}[t]{0.34\textwidth}
\vspace{0pt}
\centering
\captionsetup{justification=centerlast}
\captionof{table}{Example 1, Test 2: \mbox{Observed} values of $k$.}
\begin{tabular}{|l|c|}\hline
\multicolumn{1}{|c|}{\textbf{Method}} & $k$ \\\hline
gPC                     &11.3\\\hline
Interpolation B-spline  &5.6\\\hline
Approximation B-spline  &0.3\\\hline
SP spline               &0.6\\\hline
CWENO interpolation     &6.0\\\hline
\end{tabular}
\label{tab:3.1.2}
\end{minipage}
\end{figure}

In general, similar results are obtained as in Test 1. However, since the interval for $\xi$ is wider, more points are required to avoid abrupt jumps in CWENO interpolation. This can be observed in \fref{fig:3.1.4} where we plot the L$^1$ error for different $N$: while the convergence rates are similar to Test~1, the error magnitude for the CWENO interpolation has increased. The observed values of $k$ are summarized in~\tref{tab:3.1.2}.

\subsection{Example 2: A Discontinuous Function}
In this example we introduce discontinuity in~\eref{eq:3.1.1}:
\begin{equation*}
U(\xi)=
\begin{cases}
-3\cos(\pi\xi) & \text{if }x<0.1\\
\phantom{-}3\cos(\pi\xi) & \text{otherwise}
\end{cases}
\end{equation*}
where $\xi\sim\mathcal{U}[-1,1]$. The objective remains the same: computation of $p(U)$ and estimates for $\mu$ and $\sigma$.

\begin{figure}[ht!]
\centering 
\begin{subfigure}[t]{0.30\textwidth} 
\centering
\includegraphics[width=\textwidth]{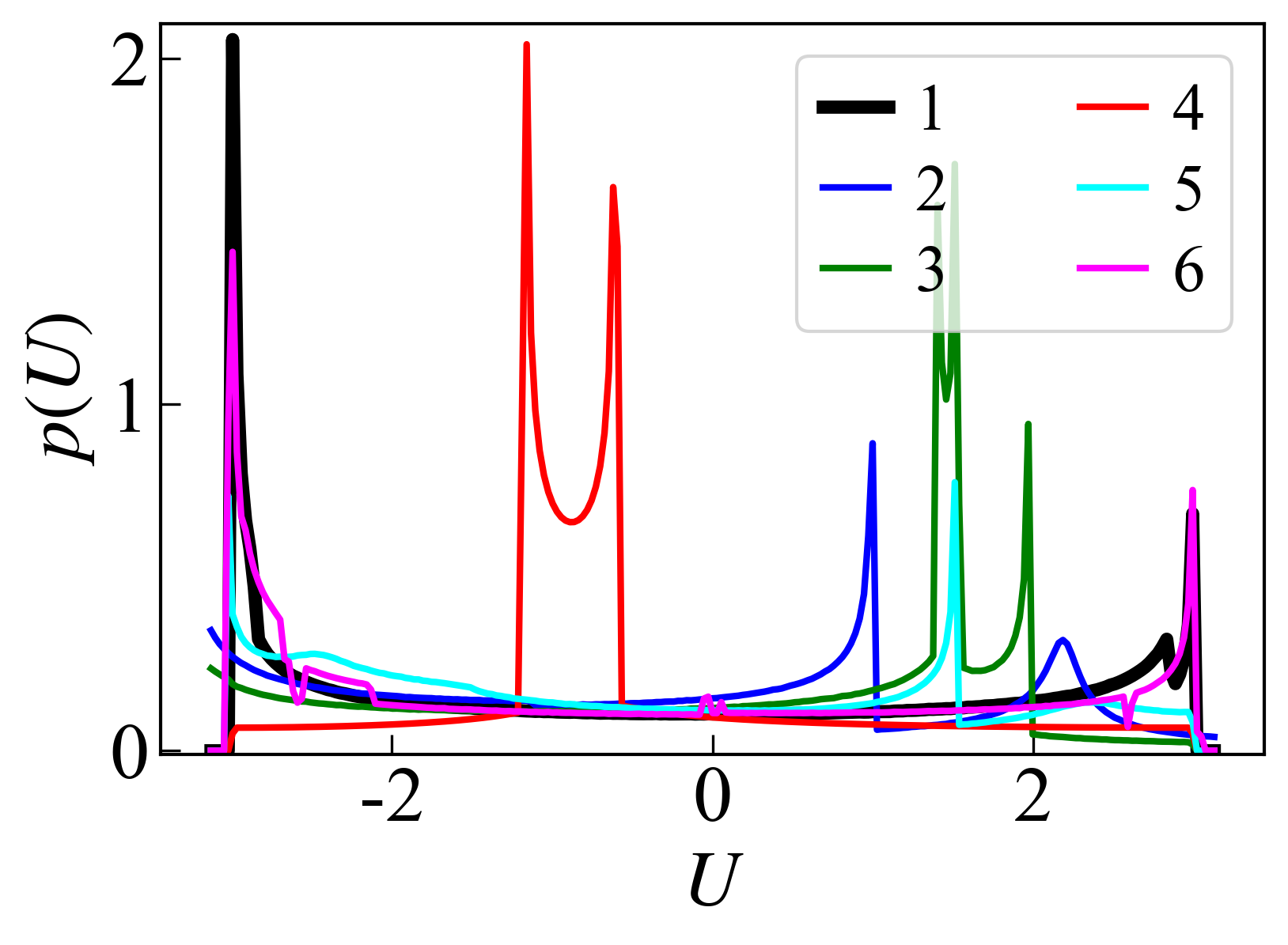}
\caption{}
\label{fig:3.2.1a}
\end{subfigure}
\begin{subfigure}[t]{0.30\textwidth}
\centering 
\includegraphics[width=\textwidth]{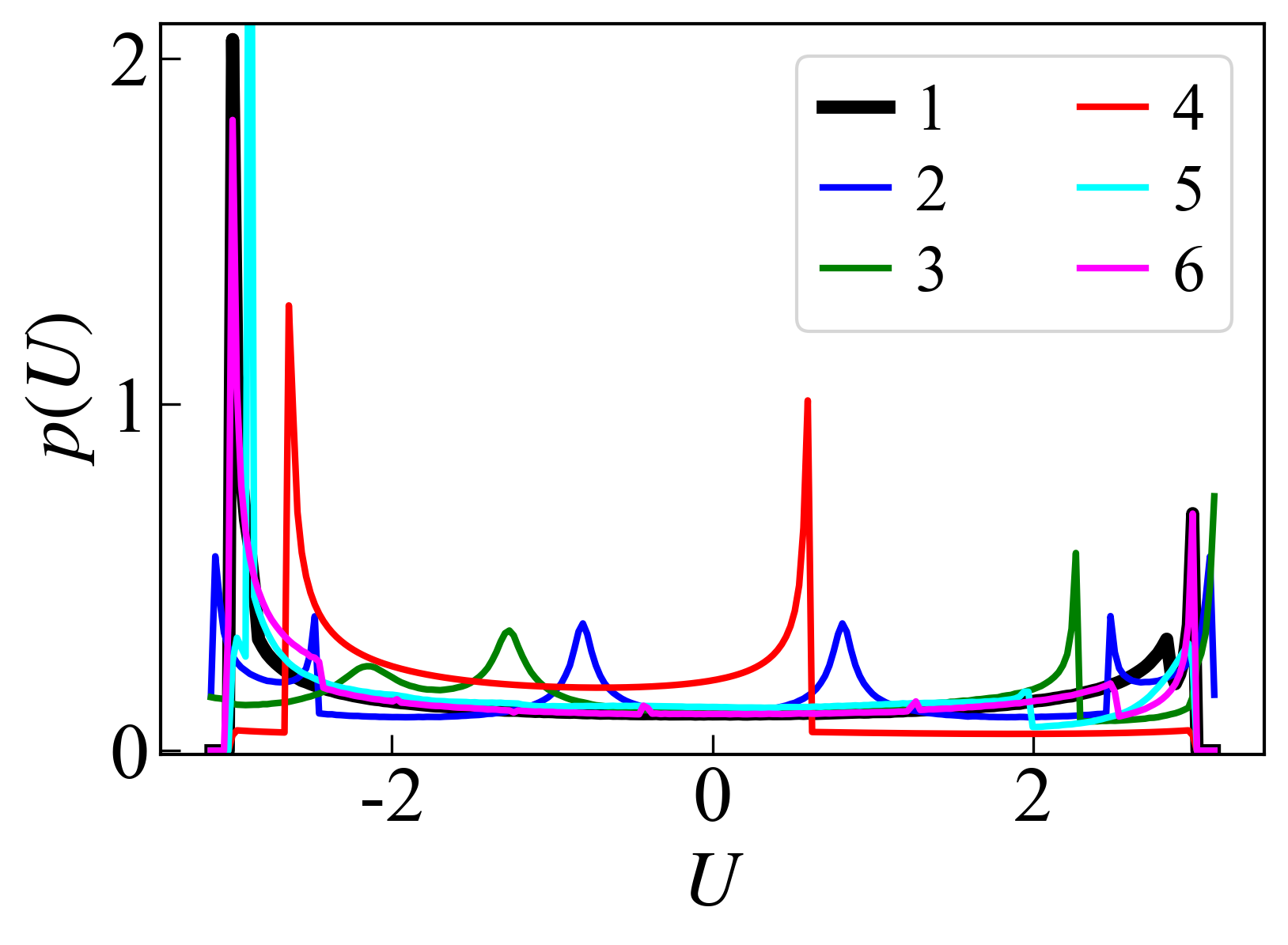}
\caption{}
\label{fig:3.2.1b}
\end{subfigure}
\begin{subfigure}[t]{0.30\textwidth}
\centering
\includegraphics[width=\textwidth]{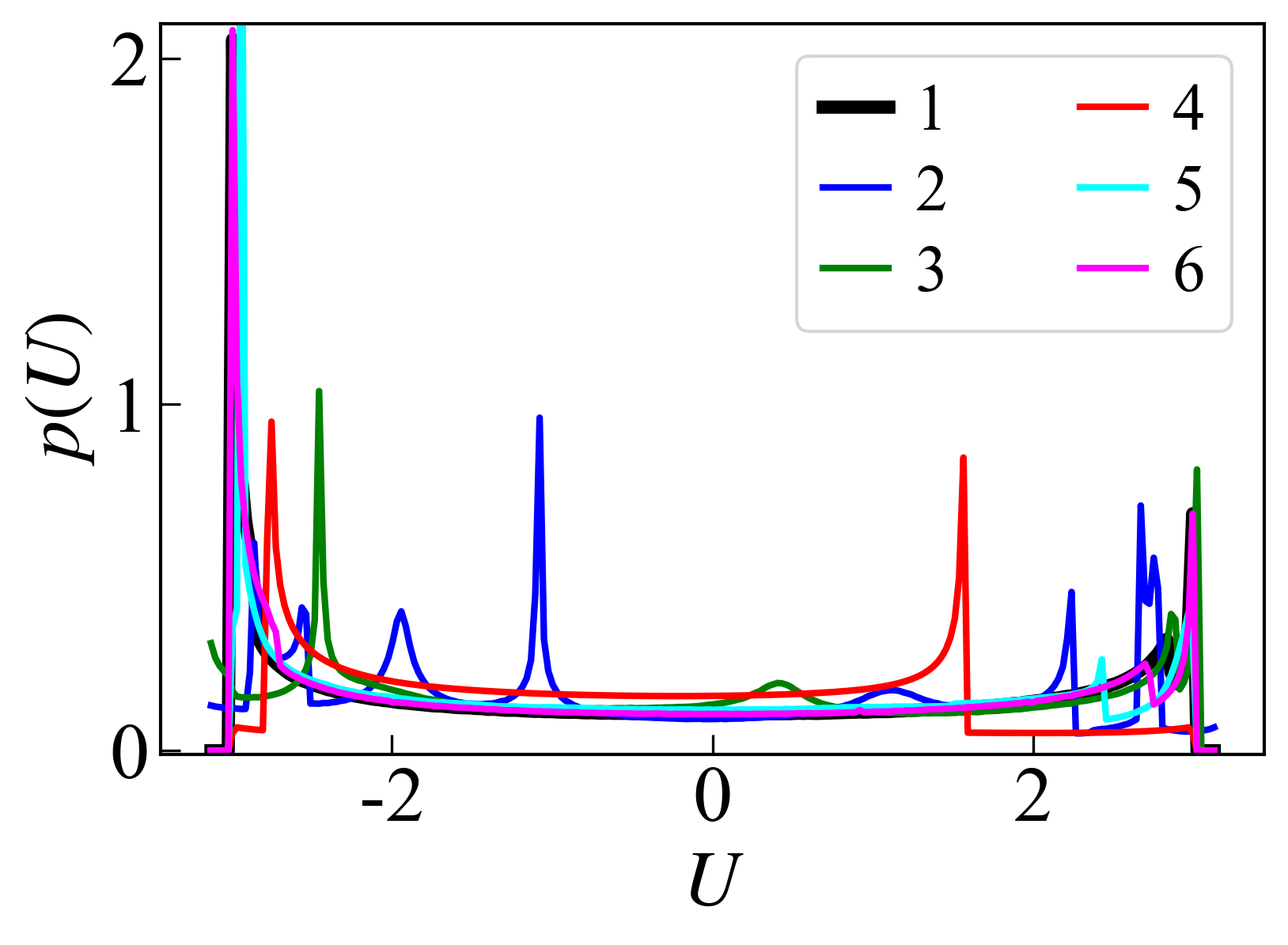}
\caption{} 
\label{fig:3.2.1c} 
\end{subfigure} 
\captionsetup{justification=centerlast} 
\caption{Example 2: PDF of $U$ for (\subref{fig:3.2.1a}) $N=7$, (\subref{fig:3.2.1b}) $N=12$, and (\subref{fig:3.2.1c}) $N=16$.} 
\label{fig:3.2.1} 
\end{figure}

We follow the same procedure as in Example~1 and generate $M$ random samples to obtain the reference PDF $\tilde{p}(U)$. Then the approximation/interpolation methods are used to reconstruct the correponding PDFs $p(U)$. \fref{fig:3.2.1} shows the resulting PDFs for $N=7$, $12$, and $16$. One can see that since for a discontinuous solution gPC and interpolation B-spline yield oscillatory behavior, their application results in the failure of restoring the PDF. As in Example~1, the approximation B-spline smears data and fails to correctly reproduce the PDF, while CWENO interpolation still requires a sufficient number of the collocation points. Although the SP spline provides satisfactory results, the CWENO interpolation outperforms it in accuracy.

\begin{figure}[ht!]
\centering
\begin{minipage}[t]{0.64\textwidth}
\vspace{0pt}
\centering
\includegraphics[width=0.5\textwidth]{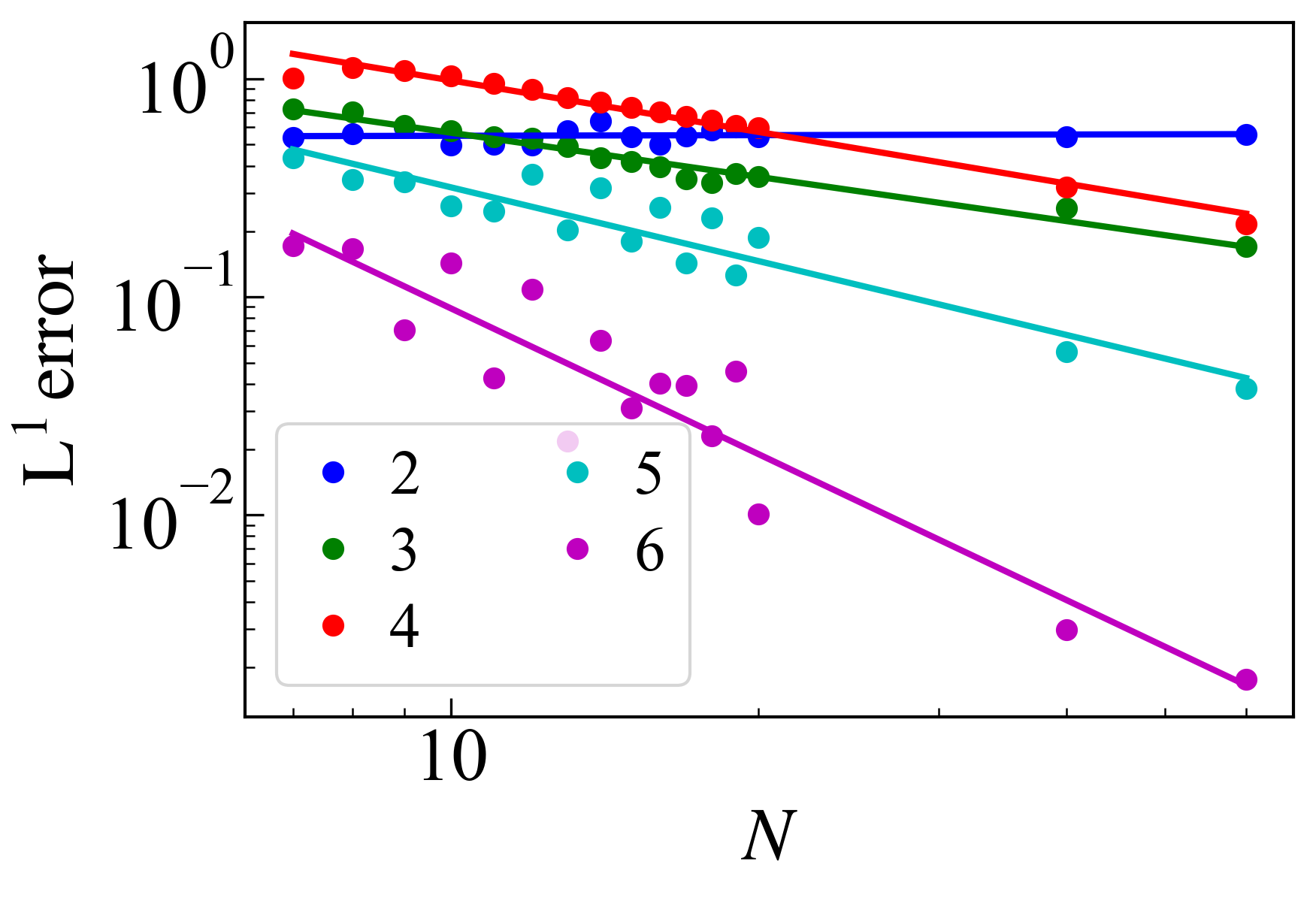}
\captionsetup{justification=centerlast}
\caption{Example 2: L$^1$ error of PDF reconstruction as a function of $N$, and the corresponding power-law fit (solid lines).}
\label{fig:3.2.2}
\end{minipage}
\hfill
\begin{minipage}[t]{0.34\textwidth}
\vspace{0pt}
\centering
\captionsetup{justification=centerlast}
\captionof{table}{Example 2: \mbox{Observed} values of $k$.}
\begin{tabular}{|l|c|}\hline
\multicolumn{1}{|c|}{\textbf{Method}} & $k$ \\\hline
gPC                     & $-$0.01 \\\hline
Interpolation B-spline  & 0.7    \\\hline
Approximation B-spline  & 0.8    \\\hline
SP spline               & 1.1    \\\hline
CWENO interpolation     & 2.2    \\\hline
\end{tabular}
\label{tab:3.2.1}
\end{minipage}
\end{figure}

\begin{figure}[ht!]
\centering 
\begin{subfigure}[t]{0.30\textwidth}
\centering
\includegraphics[width=\textwidth]{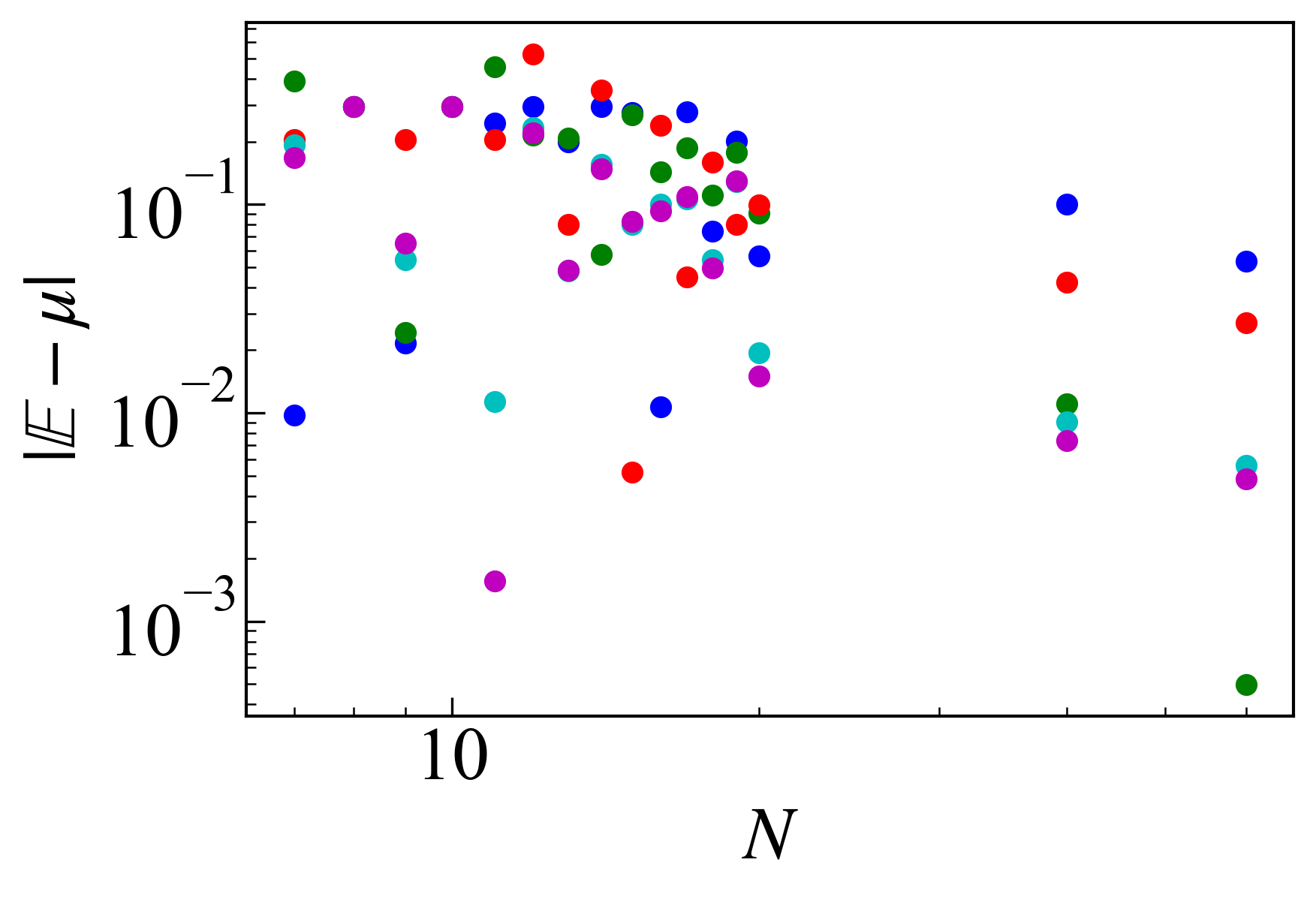}
\caption{}
\label{fig:3.2.3a}
\end{subfigure}
\begin{subfigure}[t]{0.30\textwidth}
\centering 
\includegraphics[width=\textwidth]{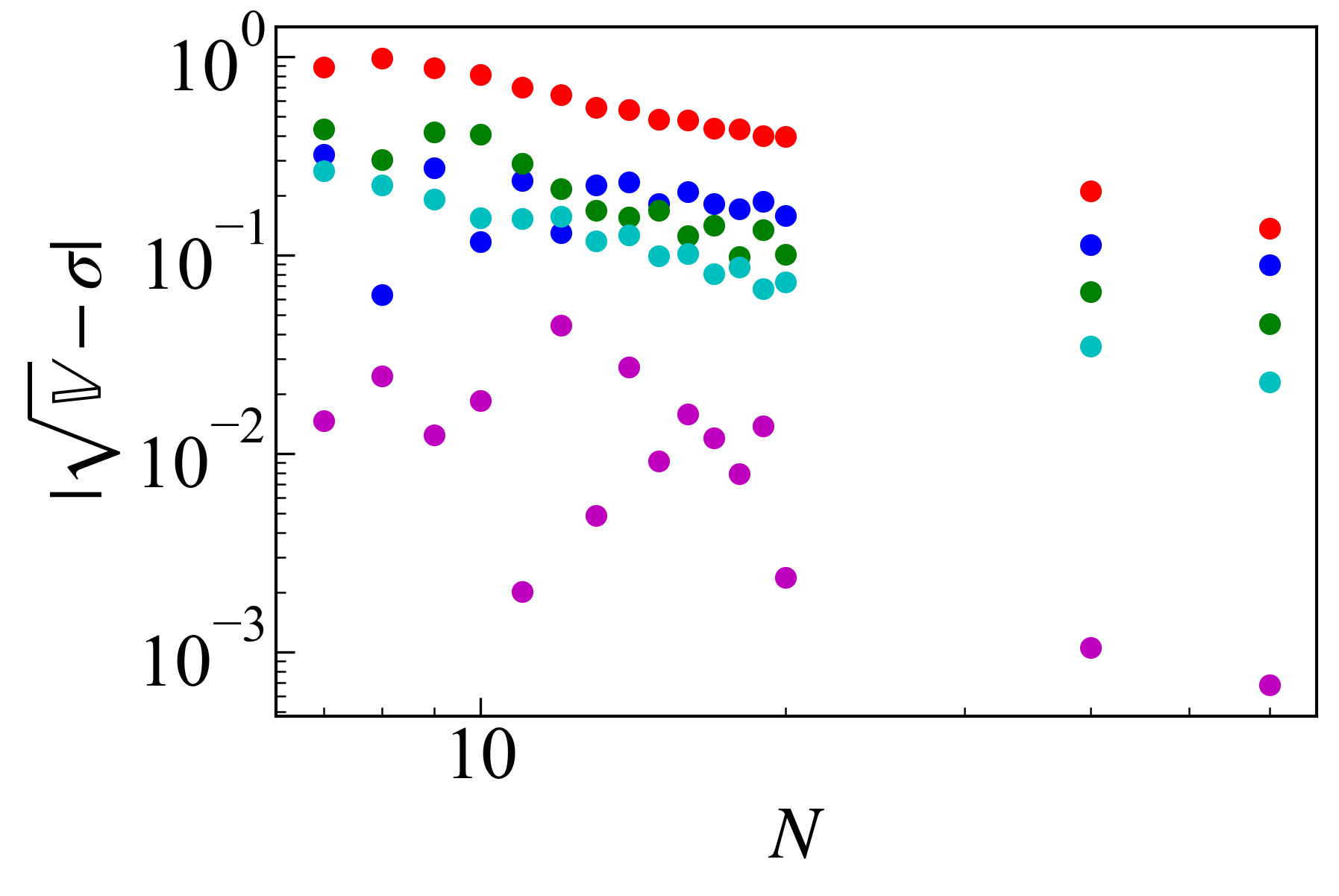}
\caption{}
\label{fig:3.2.3b}
\end{subfigure}
\captionsetup{justification=centerlast} 
\caption{Example 2: Error of (\subref{fig:3.2.3a}) mean and (\subref{fig:3.2.3b}) standard deviation estimates for different $N$.}
\label{fig:3.2.3} 
\end{figure}

\fref{fig:3.2.2} shows the L$^1$ error as a function of $N$. All methods except gPC show the reduction in error as $N$ increases and CWENO interpolation exhibits the fastest convergence rate. The observed values of $k$ are summarized in~\tref{tab:3.2.1}. \fref{fig:3.2.3} shows the absolute error of the mean and standard deviation estimates as functions of $N$. The results highlight the advantage of CWENO interpolation for discontinuous solutions. Note that while the gPC does not converge in terms of the PDF, the mean and standard deviation estimates converge.

\subsection{Example 3: Euler Equations}
In this example we consider the 1-D Euler equations of gas dynamics: 
\begin{equation}
\bm U_t+\bm F(\bm U)_x=\bm S(\bm U,x;\xi),
\label{eq:3.3.1}
\end{equation}
\begin{equation*}
\bm U=(\rho,\rho u,E)^\top,\quad\bm F(\bm U)=(\rho u,\rho u^2+P,u(E+P))^\top,\quad\bm S\equiv\bm0.
\end{equation*}
Here, $\rho(t,x;\xi)$ is the density, $u(t,x;\xi)$ is the velocity, $E(t,x;\xi)$ is the total energy, and $P(t,x;\xi)$ is the pressure related to the conservative quantities through the equation of state for the ideal gas:
\begin{equation*}
P=(\gamma-1)\left[E-\frac{\rho u^2}{2}\right],
\end{equation*}
where $\gamma$ is the specific heat ratio.

We use semi-discrete second-order central-upwind schemes from~\cite{kurganov2001} to solve the system~\eref{eq:3.3.1} $N$ times for each collocation point $\{\xi_n\}_{n=1}^{N}$. Time integration is performed using the third-order Runge-Kutta method~\cite{gottlieb2011} with the time step chosen according to the CFL number $0.45$. As this solution requires a noticeable computation time, we refrain from running $M$ simulations to restore the reference PDF.

We consider the Sod shock tube problem with $\gamma=1.4$ and initial conditions
\begin{equation}
(\rho(0,x),u(0,x),P(0,x))=\left\{\begin{aligned}
&(1,0,1),&&x\leq0.5,\\
&(0.125,0,0.1),&&x>0.5,
\end{aligned}\right.
\label{eq:3.3.4}
\end{equation}
and free BCs in the spatial domain $x\in[0,1]$. We then perturb the initial density in~\eref{eq:3.3.4} as
\begin{equation*}
\rho(0,x;\xi)=
\left\{\begin{aligned}
&1+0.1\xi,&&x\leq0.5,\\
&0.125,&&x>0.5.
\end{aligned}\right.
\end{equation*}

The solution is computed until time $T=0.1644$ on a uniform mesh with $\dx=1/200$ and $\dxi=1/50$ (for gPC the quadrature rule is used instead). The random variable follows normal distribution $\xi\sim{\cal N}(0,1/6^2)$.

\begin{figure}[ht!]
\centering
\begin{subfigure}[t]{0.26\textwidth}
\centering
\includegraphics[width=1.0\textwidth, trim=0.2cm 1.1cm 1.0cm 2.5cm, clip]{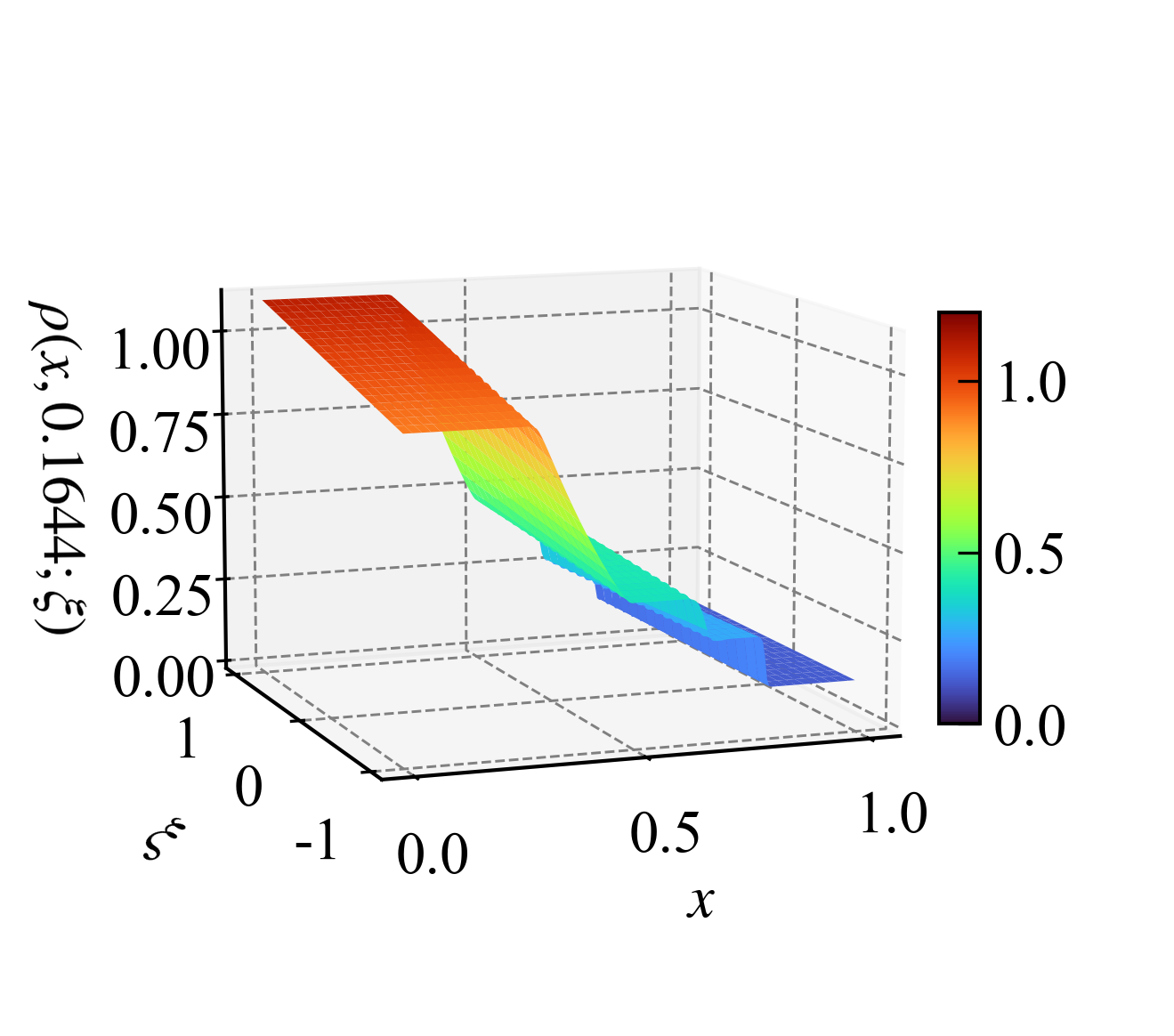}
\caption{}
\label{fig:3.3.1a}
\end{subfigure}
\hspace{0.05cm}
\centering
\begin{subfigure}[t]{0.26\textwidth}
\centering
\includegraphics[width=1.0\textwidth, trim=0.2cm 1.1cm 1.0cm 2.5cm, clip]{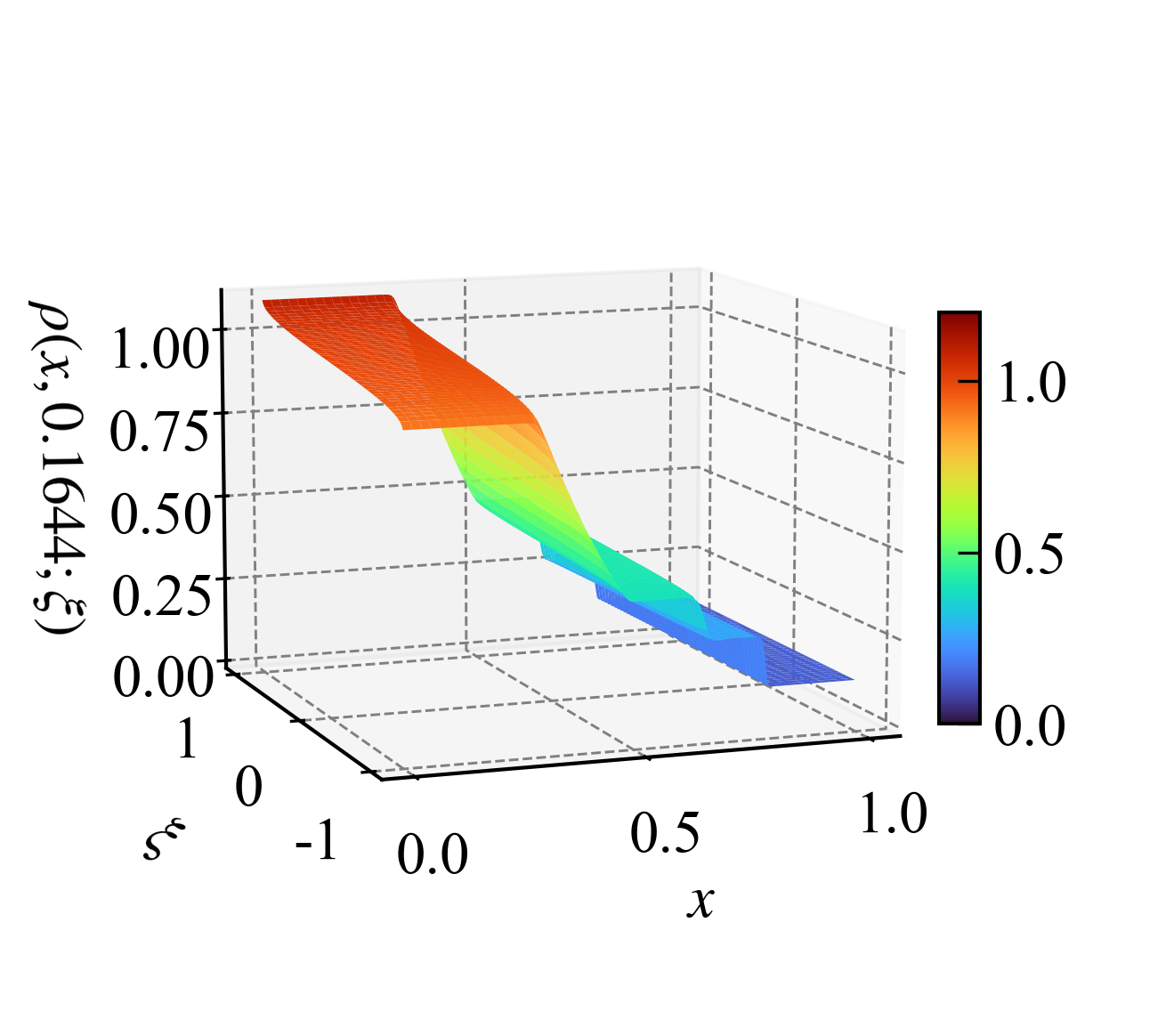}
\caption{}
\label{fig:3.3.1b}
\end{subfigure}
\hspace{0.05cm}
\centering
\begin{subfigure}[t]{0.26\textwidth}
\centering
\includegraphics[width=1.0\textwidth, trim=0.2cm 1.1cm 1.0cm 2.5cm, clip]{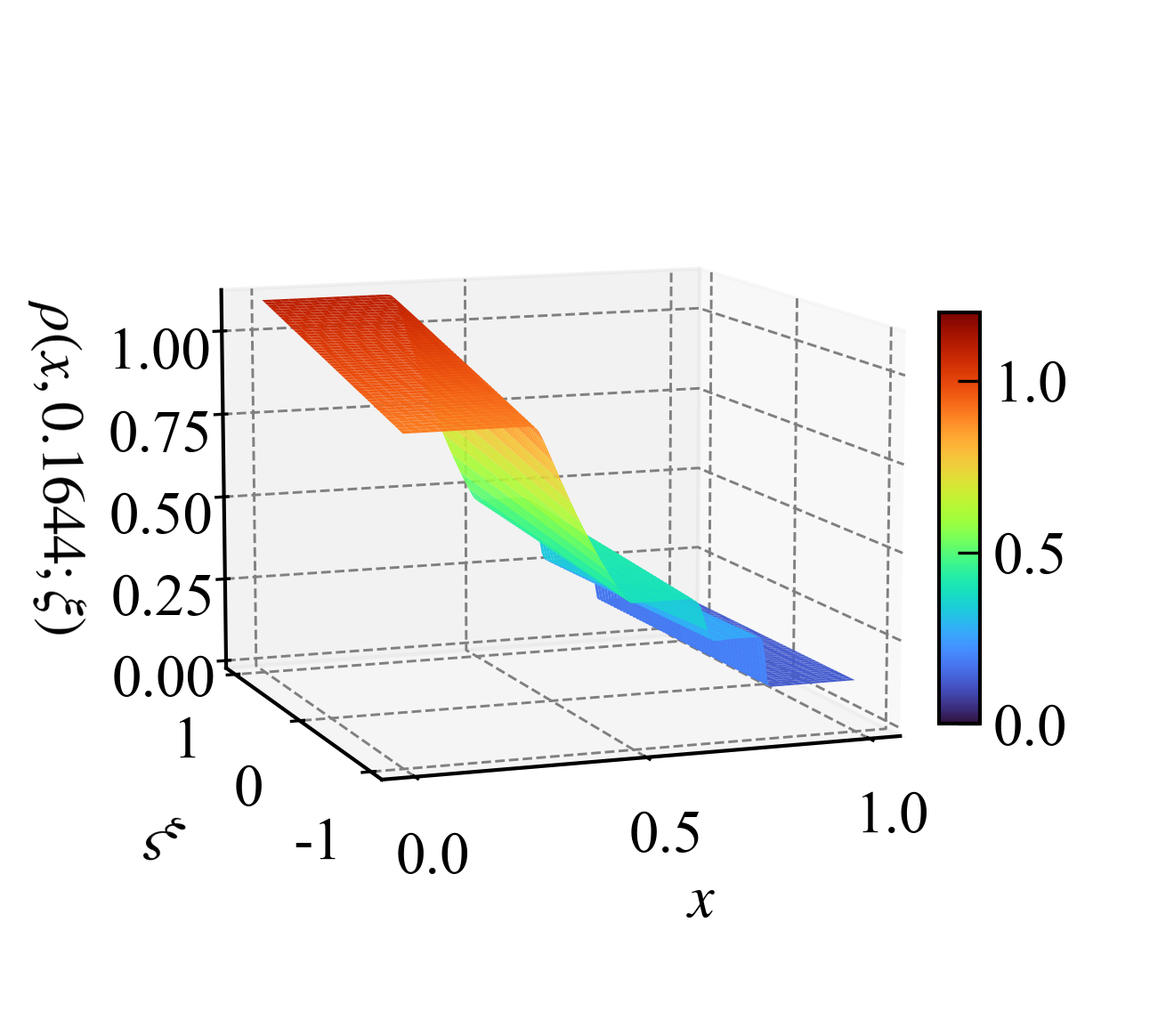}
\caption{}
\label{fig:3.3.1c}
\end{subfigure}
\centering
\begin{subfigure}[t]{0.26\textwidth}
\centering
\includegraphics[width=1.0\textwidth, trim=0.2cm 1.1cm 1.0cm 2.5cm, clip]{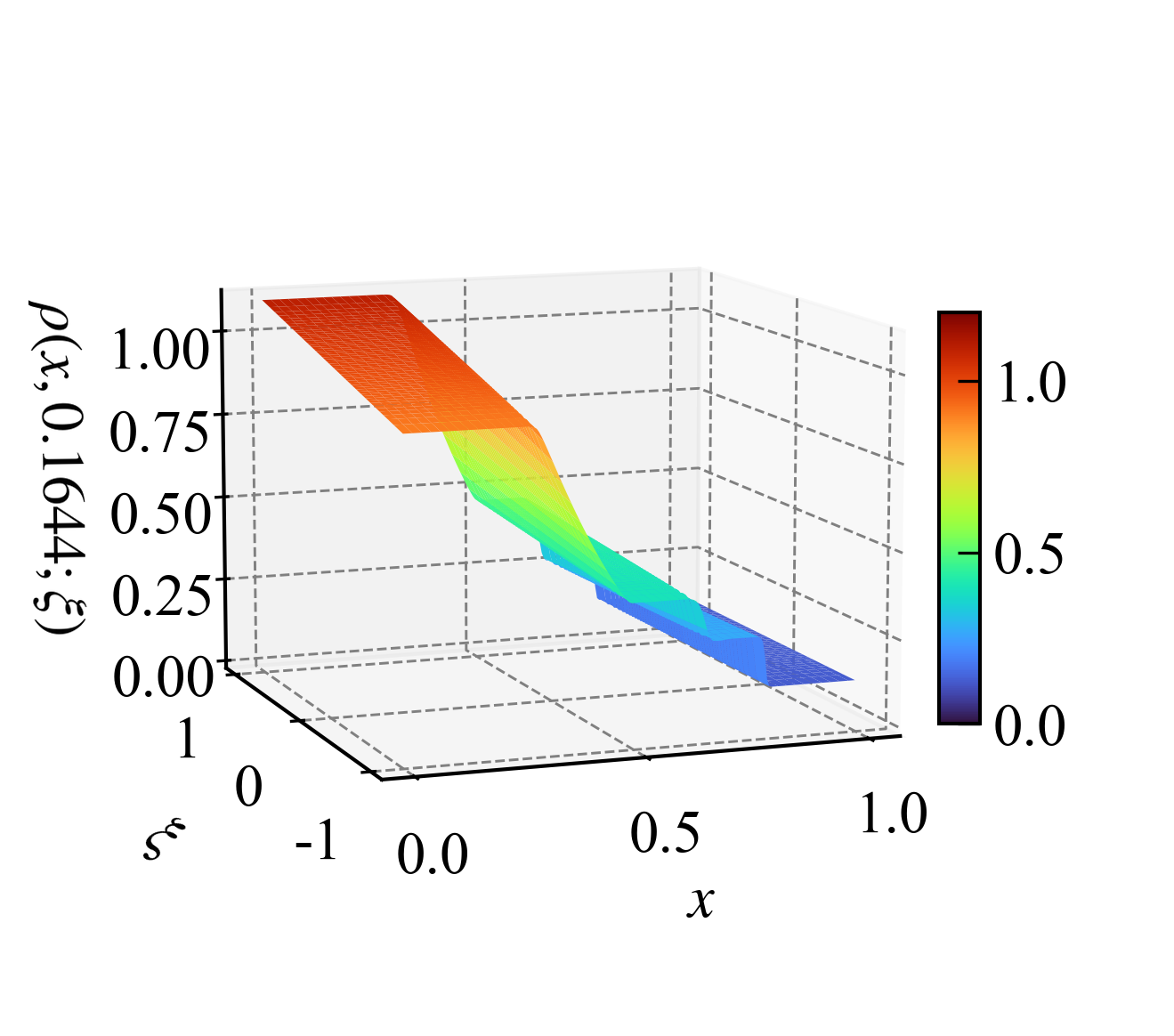}
\caption{}
\label{fig:3.3.1d}
\end{subfigure}
\hspace{0.05cm}
\centering
\begin{subfigure}[t]{0.26\textwidth}
\centering
\includegraphics[width=1.0\textwidth, trim=0.2cm 1.1cm 1.0cm 2.5cm, clip]{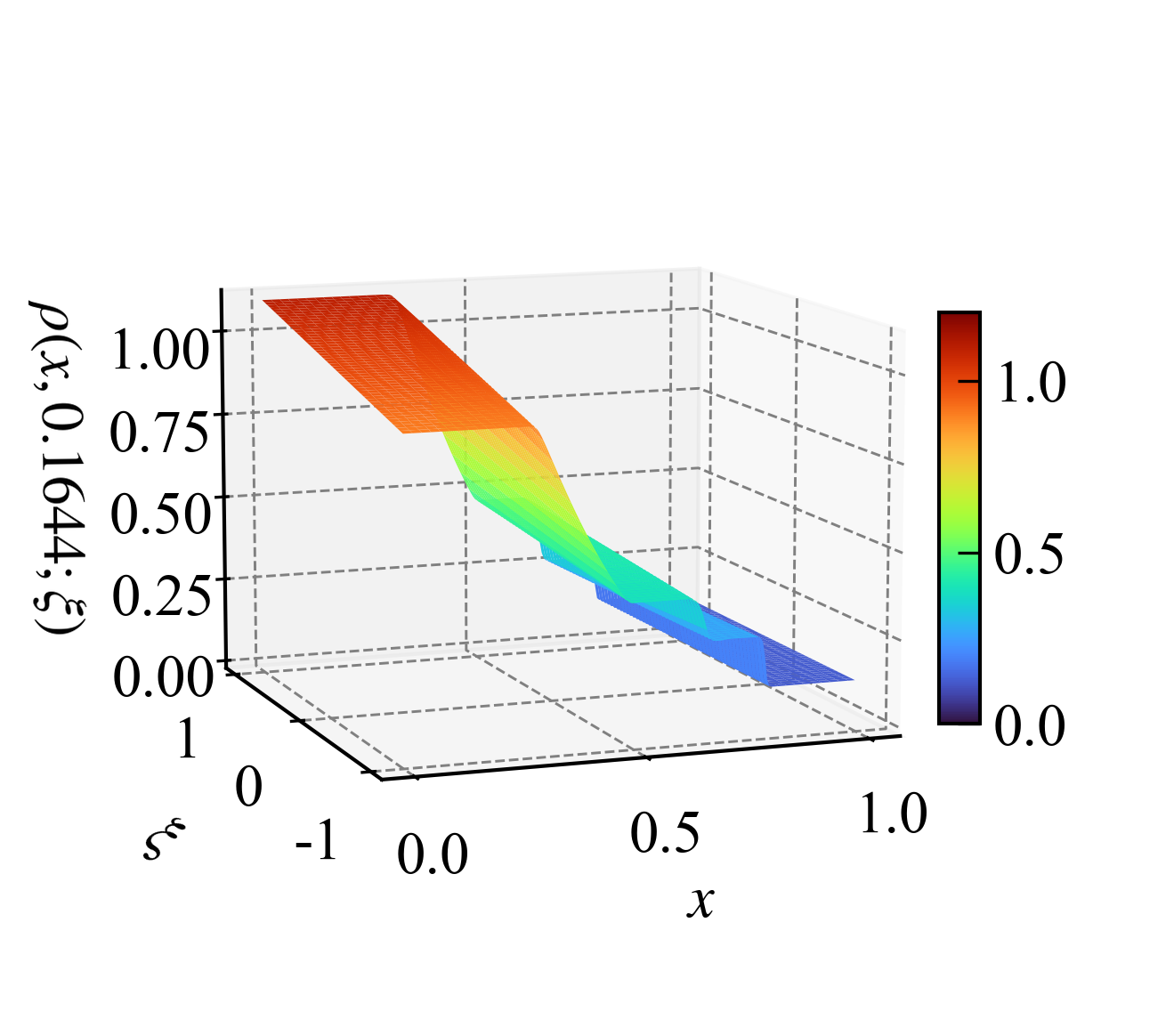}
\caption{}
\label{fig:3.3.1e}
\end{subfigure}
\hspace{0.05cm}
\centering
\begin{subfigure}[t]{0.26\textwidth}
\centering
\includegraphics[width=1.0\textwidth, trim=0.2cm 1.1cm 1.0cm 2.5cm, clip]{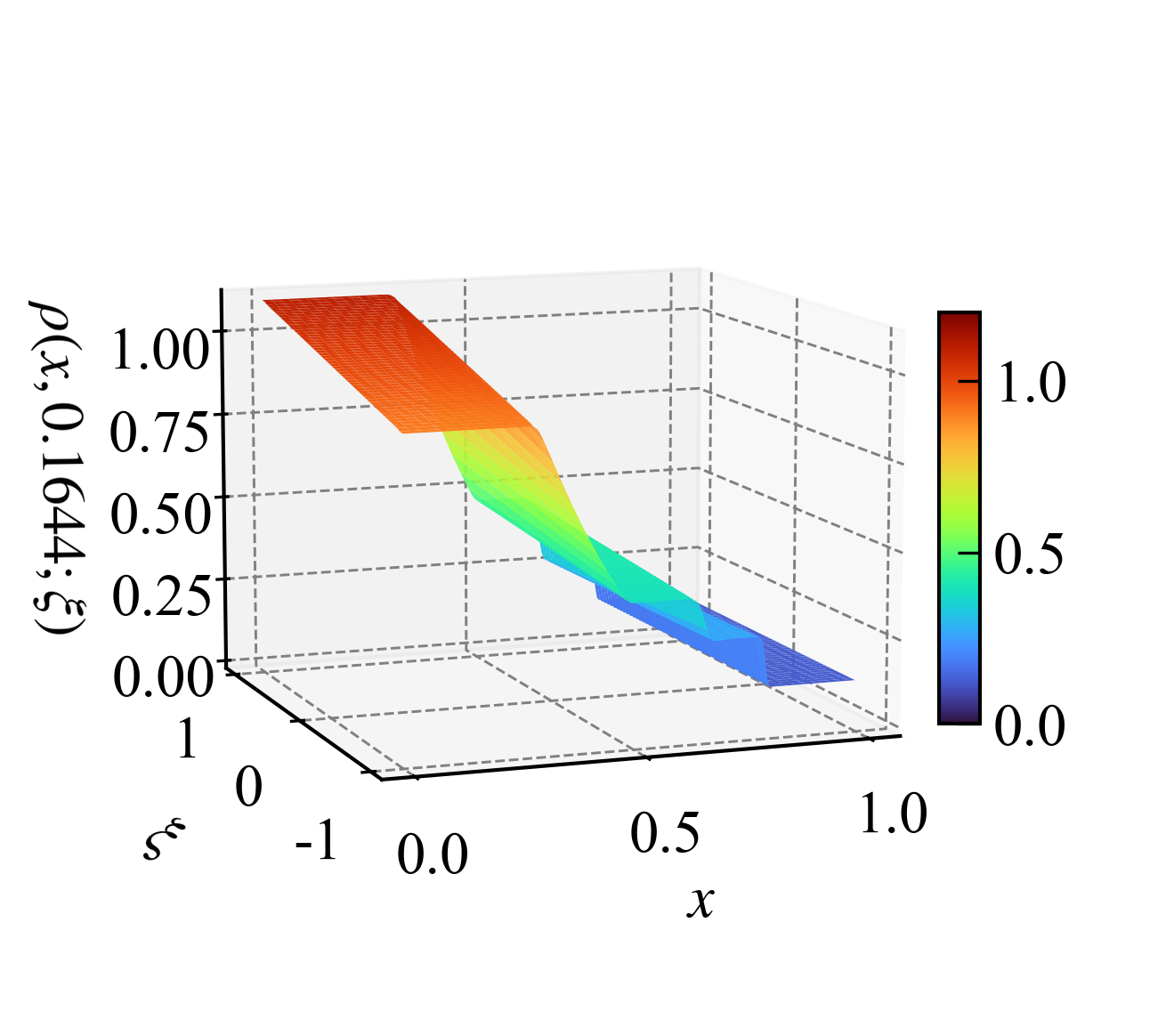}
\caption{}
\label{fig:3.3.1f}
\end{subfigure}
\captionsetup{justification=centering}
\caption{Example 3: $\rho(x,0.1644;\xi)$ obtained in (\protect\subref{fig:3.3.1a})~simulations, (\protect\subref{fig:3.3.1b})~gPC, (\protect\subref{fig:3.3.1c})~interpolation B-spline, (\protect\subref{fig:3.3.1d})~approximation B-spline, (\protect\subref{fig:3.3.1e})~SP spline, and (\protect\subref{fig:3.3.1f})~CWENO.}
\label{fig:3.3.1}
\end{figure}

The computed density is plotted in~\fref{fig:3.3.1a}. Based on this, we perform the interpolation/approximation using the proposed methods and plot the resulting surfaces in the rest of subfigures of~\fref{fig:3.3.1}. We can see that all the methods produce reasonable results. The only issue is in gPC when the tails of the distribution are extrapolated outside of the quadrature points in $\xi$,~\fref{fig:3.3.1b}.

\begin{figure}[ht!]
\centering
\begin{subfigure}[t]{0.26\textwidth}
\centering
\includegraphics[width=1.0\textwidth, trim=0.2cm 0.2cm 0.2cm 1.2cm, clip]{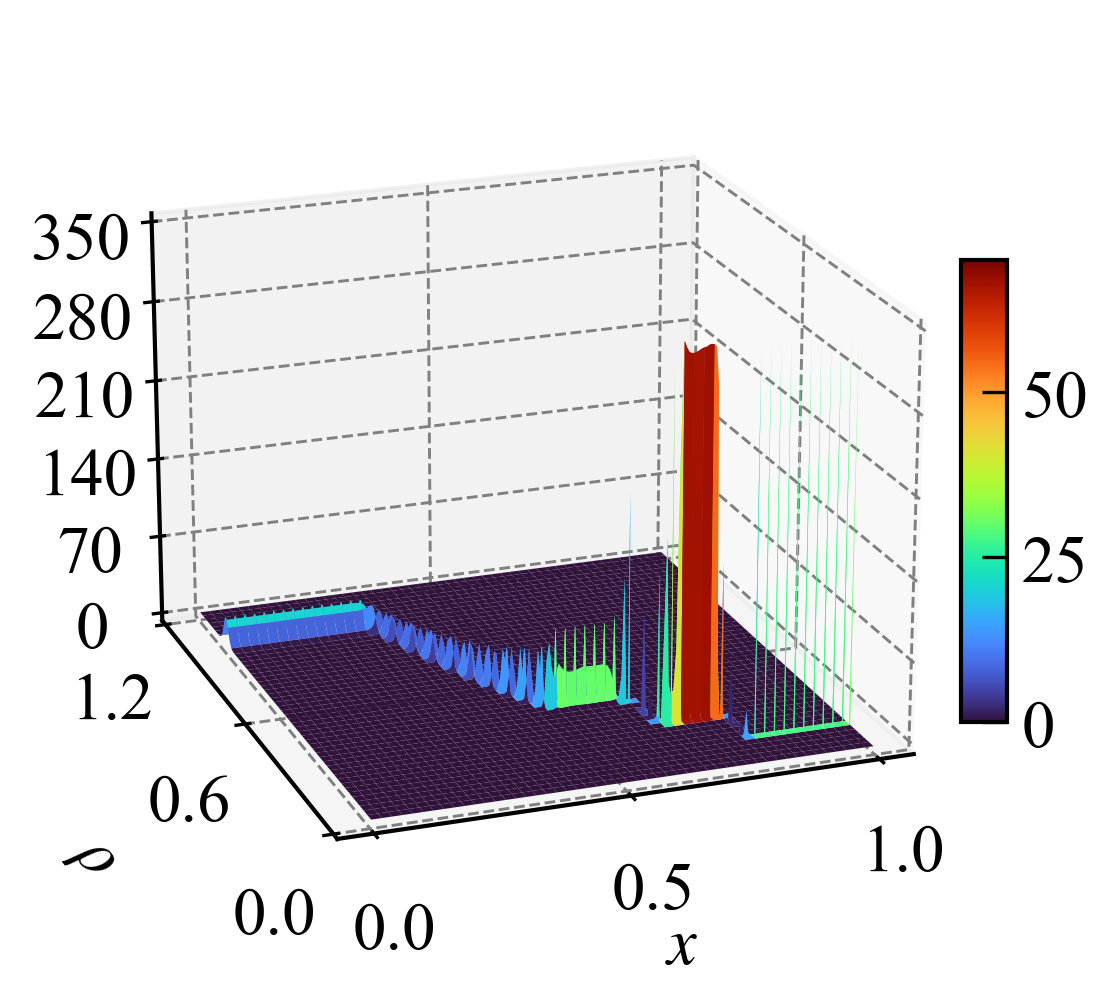}
\caption{}
\label{fig:3.3.2a}
\end{subfigure}
\centering
\begin{subfigure}[t]{0.26\textwidth}
\centering
\includegraphics[width=1.0\textwidth, trim=0.2cm 0.2cm 0.2cm 1.2cm, clip]{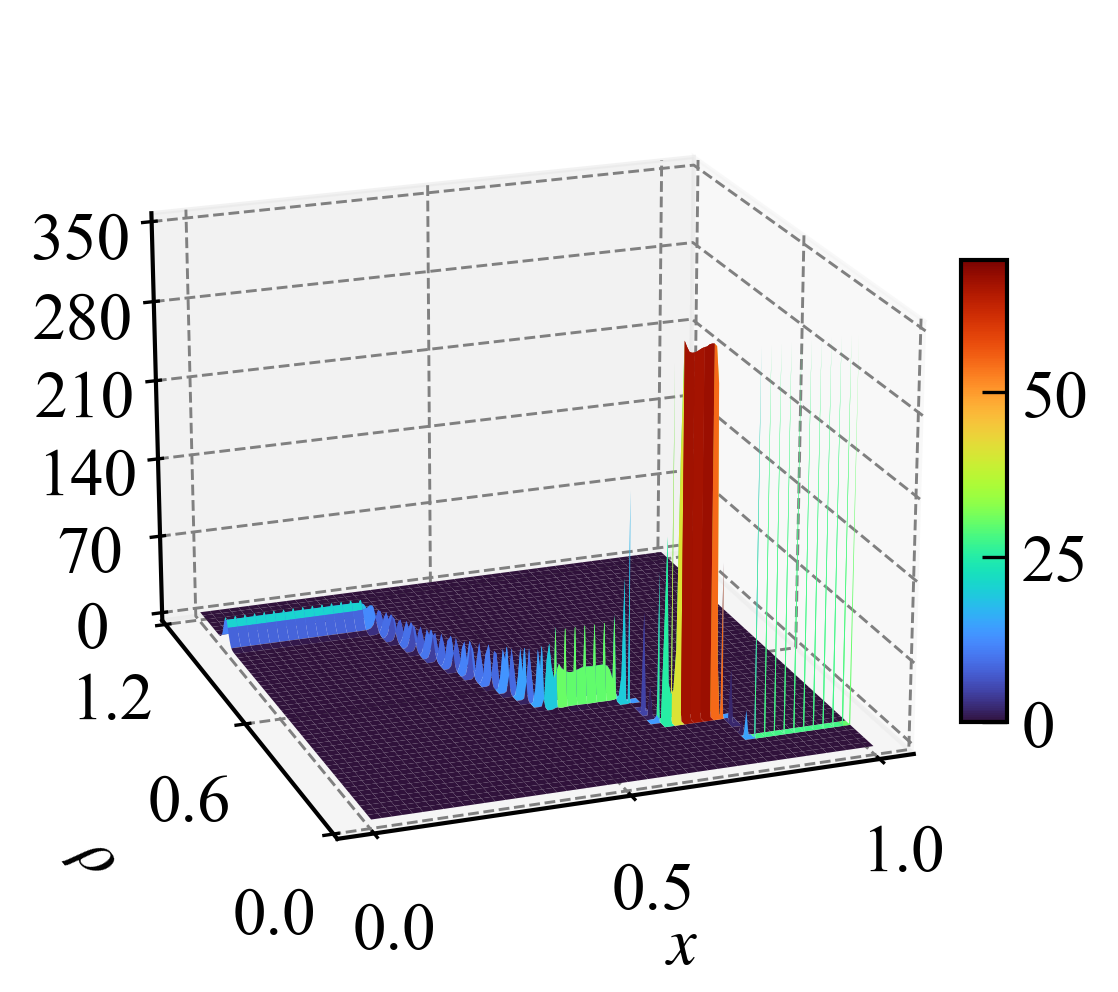}
\caption{}
\label{fig:3.3.2b}
\end{subfigure}
\centering
\begin{subfigure}[t]{0.26\textwidth}
\centering
\includegraphics[width=1.0\textwidth, trim=0.2cm 0.2cm 0.2cm 1.2cm, clip]{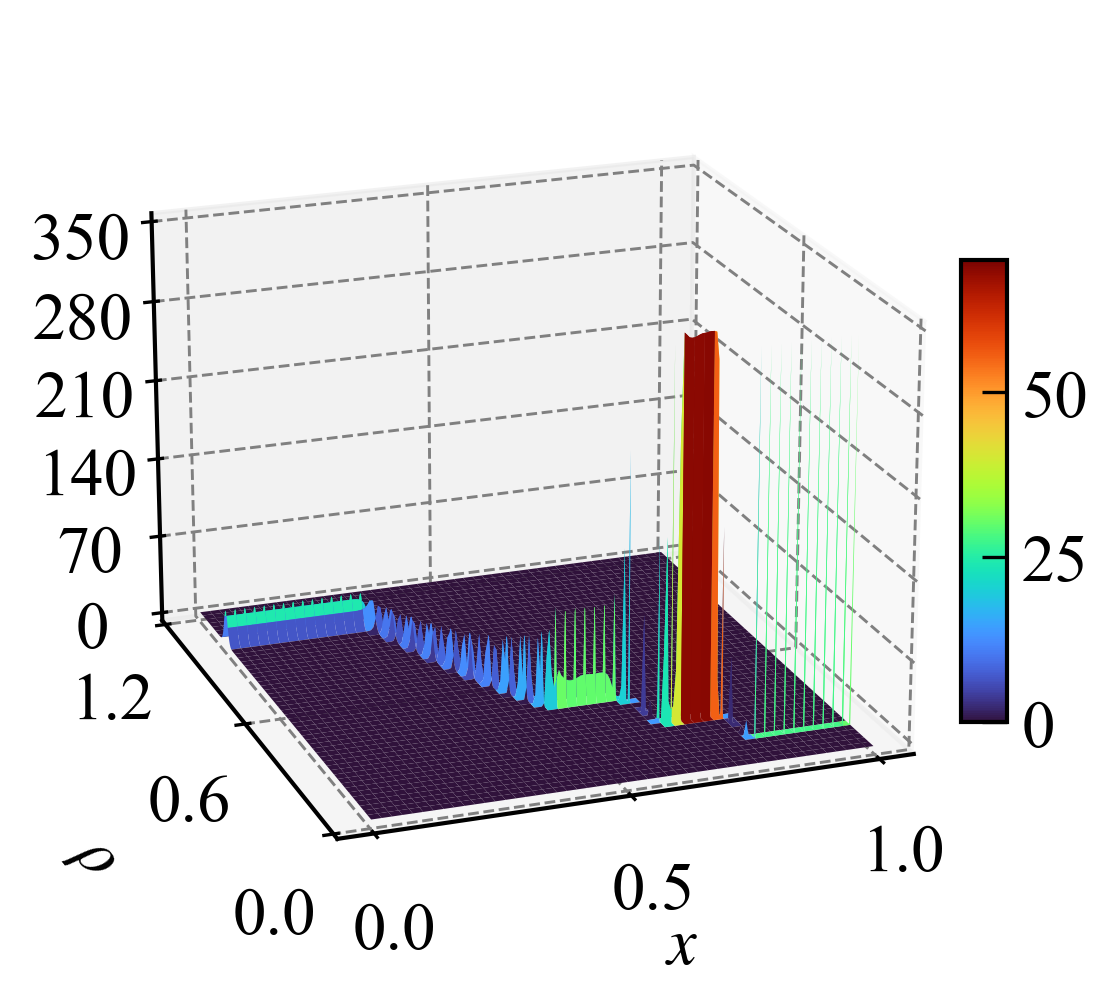}
\caption{}
\label{fig:3.3.2c}
\end{subfigure}
\newline
\centering
\begin{subfigure}[t]{0.26\textwidth}
\centering
\includegraphics[width=1.0\textwidth, trim=0.2cm 0.2cm 0.2cm 1.2cm, clip]{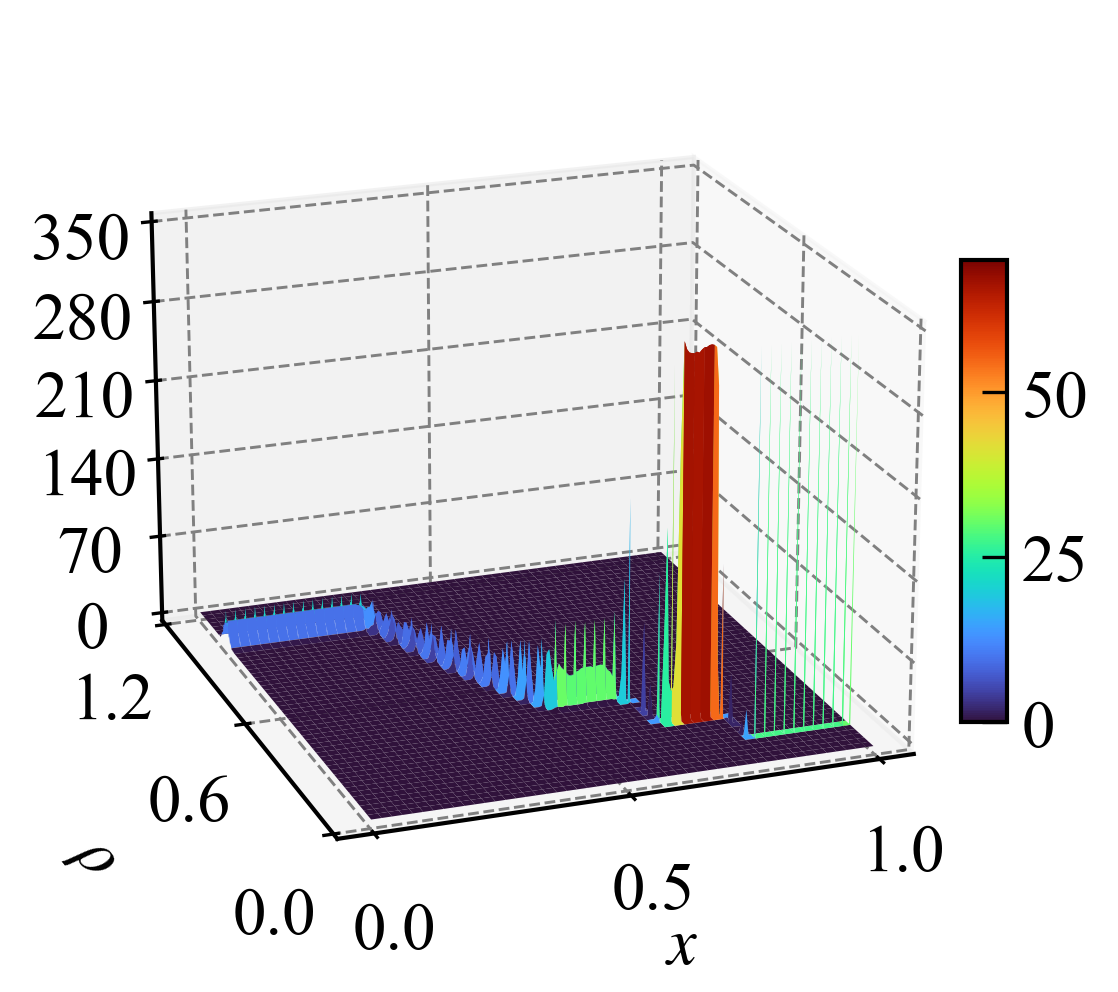}
\caption{}
\label{fig:3.3.2d}
\end{subfigure}
\centering
\begin{subfigure}[t]{0.26\textwidth}
\centering
\includegraphics[width=1.0\textwidth, trim=0.2cm 0.2cm 0.2cm 1.2cm, clip]{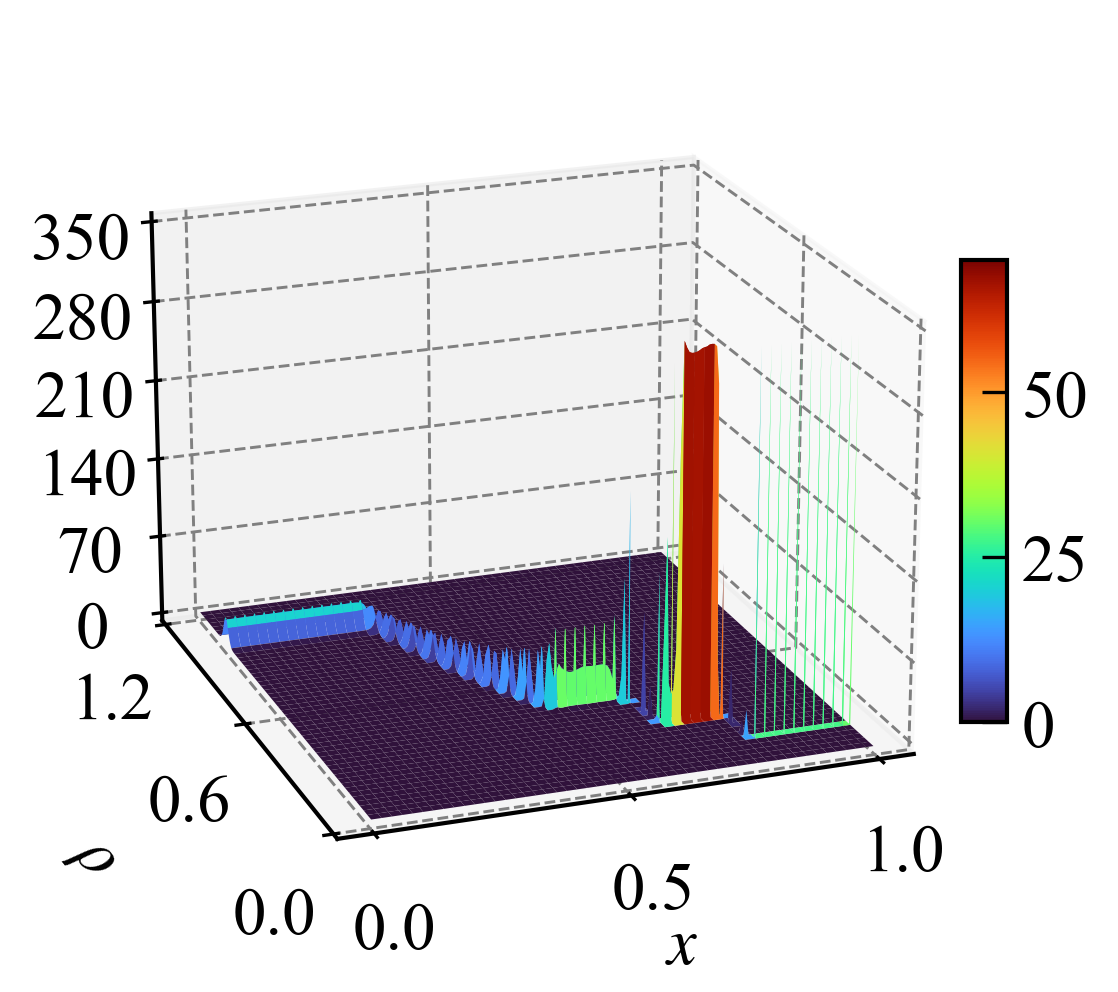}
\caption{}
\label{fig:3.3.2e}
\end{subfigure}
\captionsetup{justification=centering}
\caption{Example 3: PDF of $\rho(x,0.1644;\xi)$ obtained using (\protect\subref{fig:3.3.2a})~gPC, (\protect\subref{fig:3.3.2b})~interpolation B-spline, (\protect\subref{fig:3.3.2c})~approximation B-spline, (\protect\subref{fig:3.3.2d})~SP spline, and (\protect\subref{fig:3.3.2e})~CWENO.}
\label{fig:3.3.2}
\end{figure}

\begin{figure}[ht!]
\centering
\begin{subfigure}[t]{0.30\textwidth}
\centering
\includegraphics[width=1.0\textwidth, trim=0cm 0cm 0cm 0cm, clip]{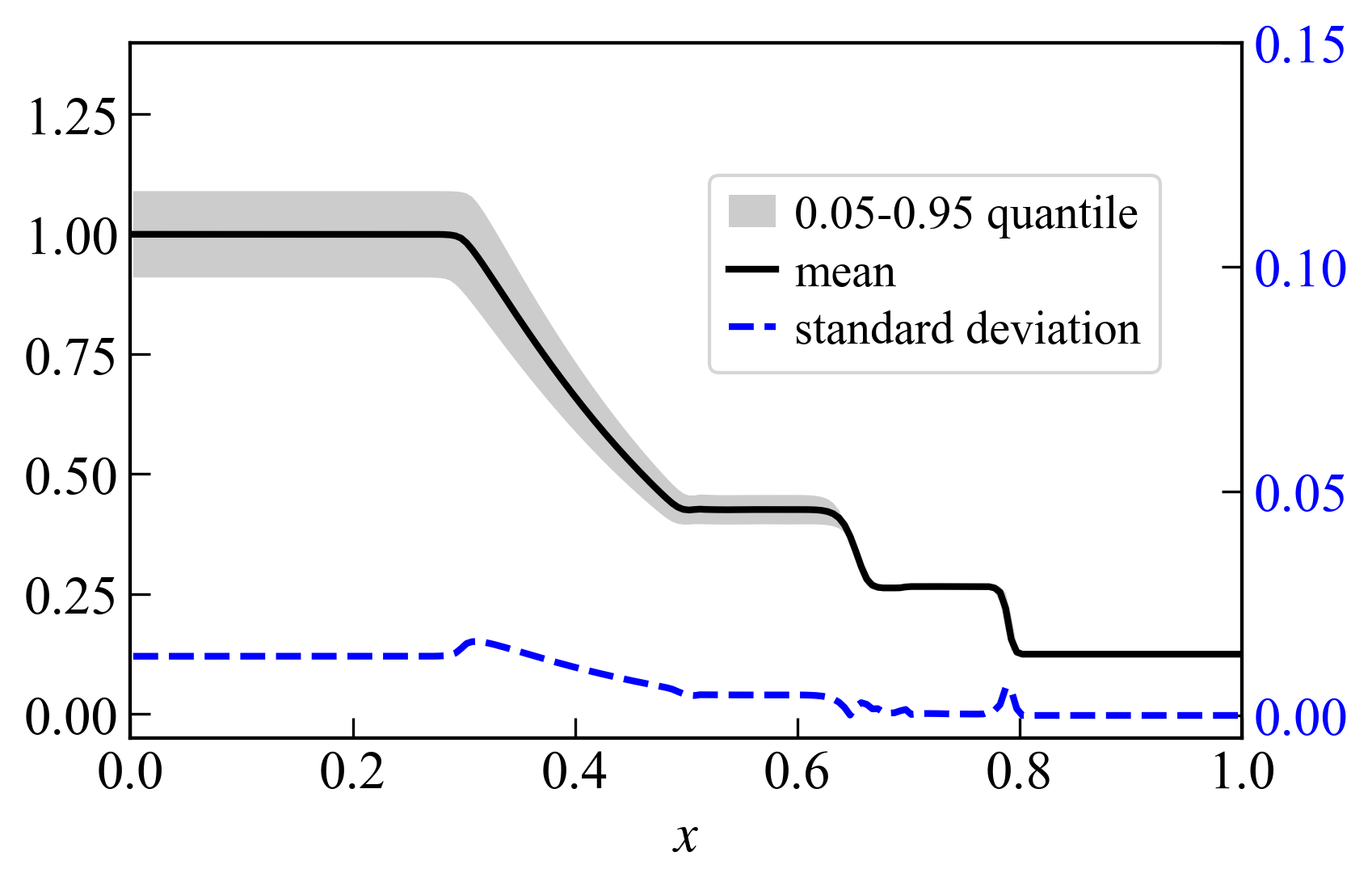}
\caption{}
\label{fig:3.3.3c}
\end{subfigure}
\begin{subfigure}[t]{0.30\textwidth}
\centering
\includegraphics[width=1.0\textwidth, trim=0cm 0cm 0cm 0cm, clip]{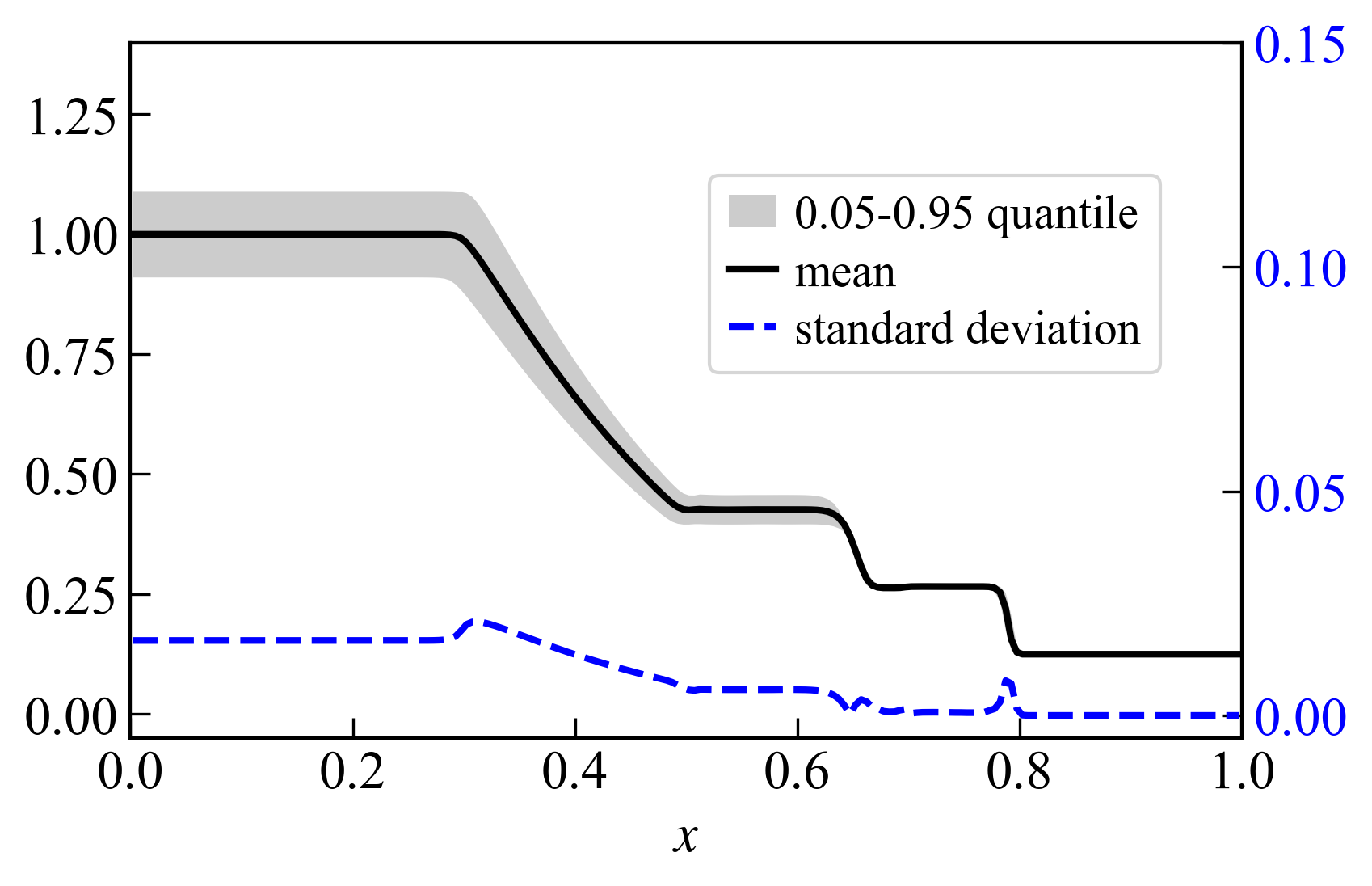}
\caption{}
\label{fig:3.3.3e}
\end{subfigure}
\captionsetup{justification=centering}
\caption{Example 3: Mean and standard deviation of $\rho(x,0.1644;\xi)$ obtained using (\protect\subref{fig:3.3.3c})~approximation B-spline and (\protect\subref{fig:3.3.3e})~CWENO.}
\label{fig:3.3.3}
\end{figure}

The interpolated data are then used to compute the PDF, $p(\rho)$, as shown in~\fref{fig:3.3.2}. All methods produce qualitatively similar results, except for the approximation B-spline, which diverges due to its smoothing properties. The PDF is subsequently used to compute the mean and standard deviation for each method, with results shown in~\fref{fig:3.3.3} for the approximation B-spline and CWENO interpolation. Notably, CWENO produces results similar to those of gPC, interpolation B-spline, and SP spline.

\subsection{Example 4: Shallow Water Equations}
In this example, we consider the Saint-Venant system of shallow water equations given by~\eref{eq:3.3.1} with
\begin{equation*}
\bm U=(h,hu)^\top,\quad\bm F(\bm U)=\Big(hu,hu^2+\frac{g}{2}h^2\Big)^\top,\quad\bm S=(0,-ghZ_x)^\top,
\end{equation*}
where $h(x,t;\xi)$ is the water depth, $u(x,t;\xi)$ is the velocity, $Z(x;\xi)$ is the bottom topography, and $g$ is the gravitational acceleration set to be $g=1$. 

The Saint-Venant system is considered in the physical domain $x\in[-1,1]$ subject to free BCs and deterministic initial data for the water surface $w=h+Z$ and velocity $u$,
\begin{equation*}
w(x,0;\xi)=\left\{\begin{aligned}&1,&&x<0,\\&0.5,&&x>0,\end{aligned}\right.\qquad u(x,0;\xi)=0,
\end{equation*}
and stochastic bottom topography
\begin{equation*}
Z(x)=\begin{cases}
0.125\xi+0.125(\cos(5\pi x)+2), &|x|<0.2,\\
0.125\xi+0.125,&\mbox{otherwise}.
\end{cases}
\end{equation*}

The system~\eref{eq:3.3.1} is numerically solved on a uniform spatial mesh consisting of cells of size $\dx=1/400$ until the final time $T=0.8$ for each of the $N=16$ collocation points. In this example $\xi\sim\mathcal{U}[-1,1]$.

\begin{figure}[ht!]
\centering
\begin{subfigure}[t]{0.26\textwidth}
\centering
\includegraphics[width=1.0\textwidth, trim=0.2cm 1.1cm 1.0cm 2.5cm, clip]{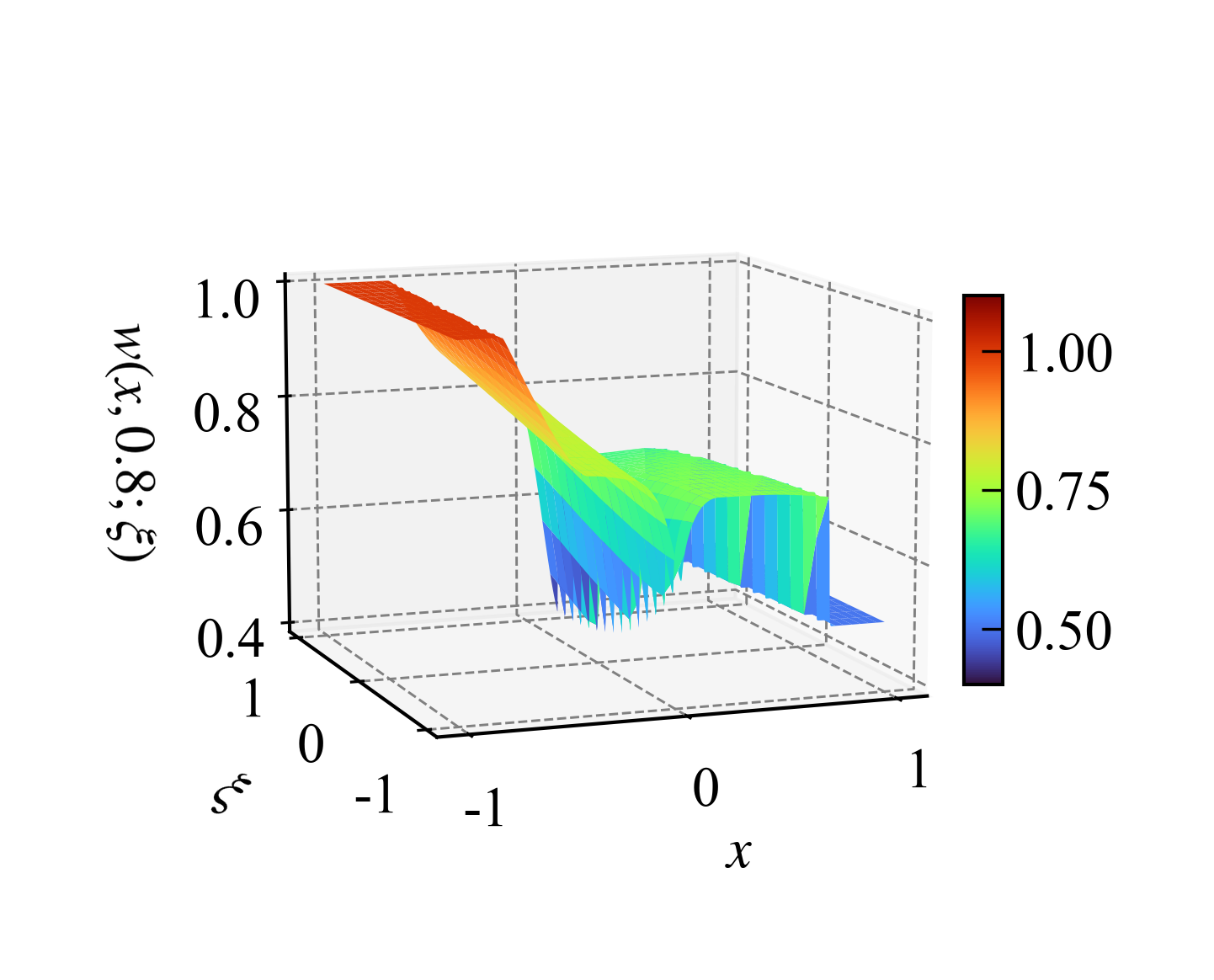}
\caption{}
\label{fig:3.4.1a}
\end{subfigure}
\hspace{0.05cm}
\centering
\begin{subfigure}[t]{0.26\textwidth}
\centering
\includegraphics[width=1.0\textwidth, trim=0.2cm 1.1cm 1.0cm 2.5cm, clip]{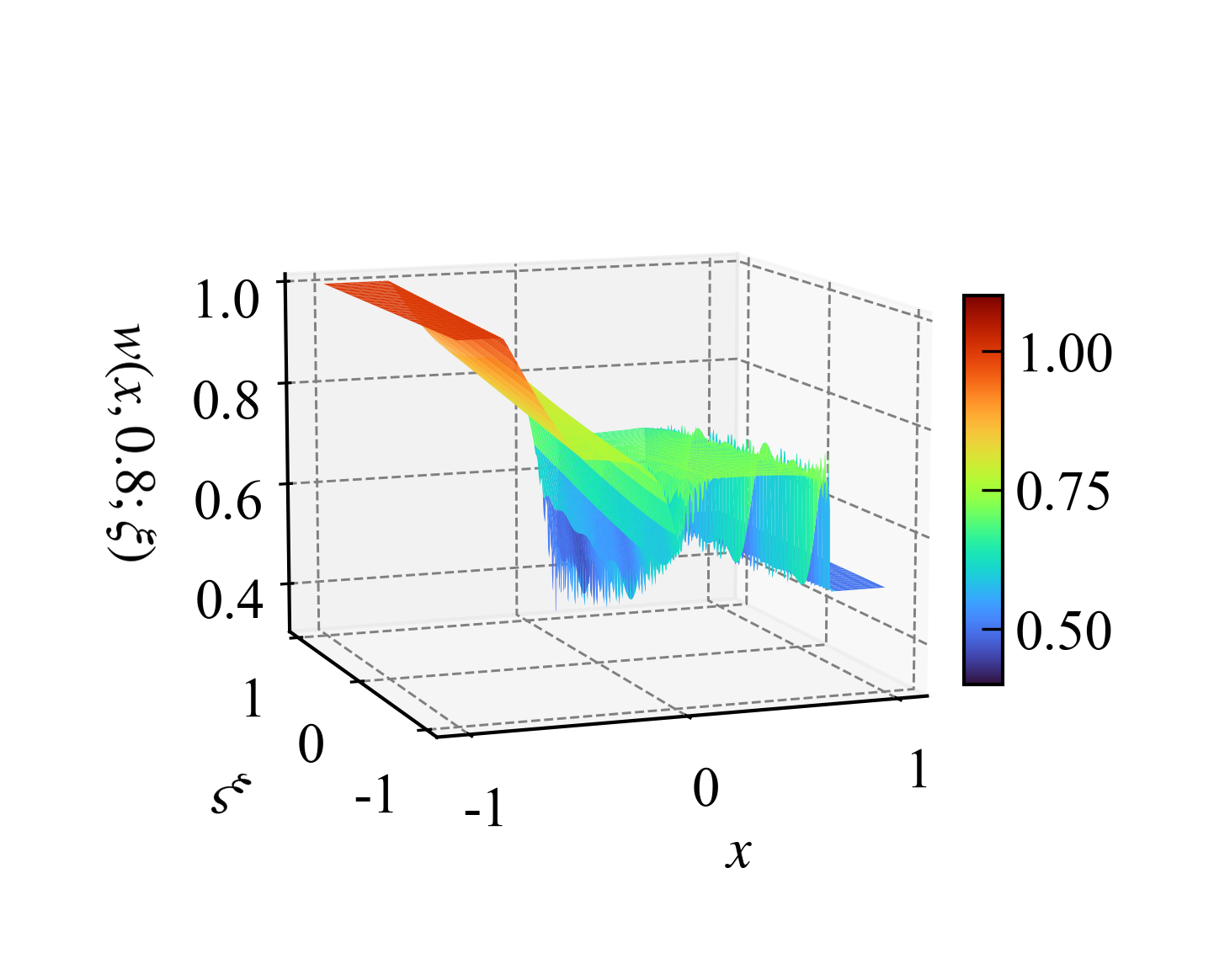}
\caption{}
\label{fig:3.4.1b}
\end{subfigure}
\hspace{0.05cm}
\centering
\begin{subfigure}[t]{0.26\textwidth}
\centering
\includegraphics[width=1.0\textwidth, trim=0.2cm 1.1cm 1.0cm 2.5cm, clip]{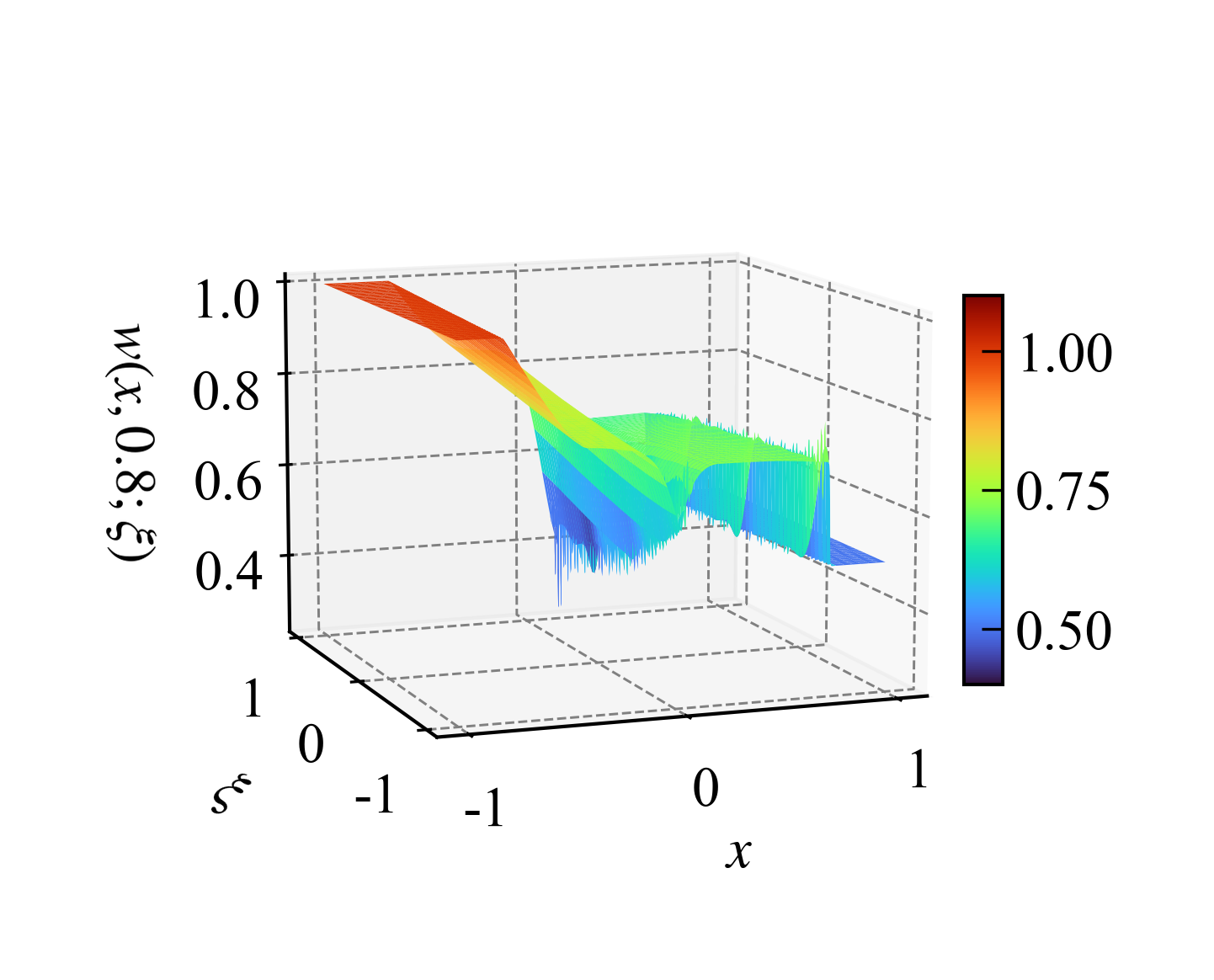}
\caption{}
\label{fig:3.4.1c}
\end{subfigure}
\centering
\begin{subfigure}[t]{0.26\textwidth}
\centering
\includegraphics[width=1.0\textwidth, trim=0.2cm 1.1cm 1.0cm 2.5cm, clip]{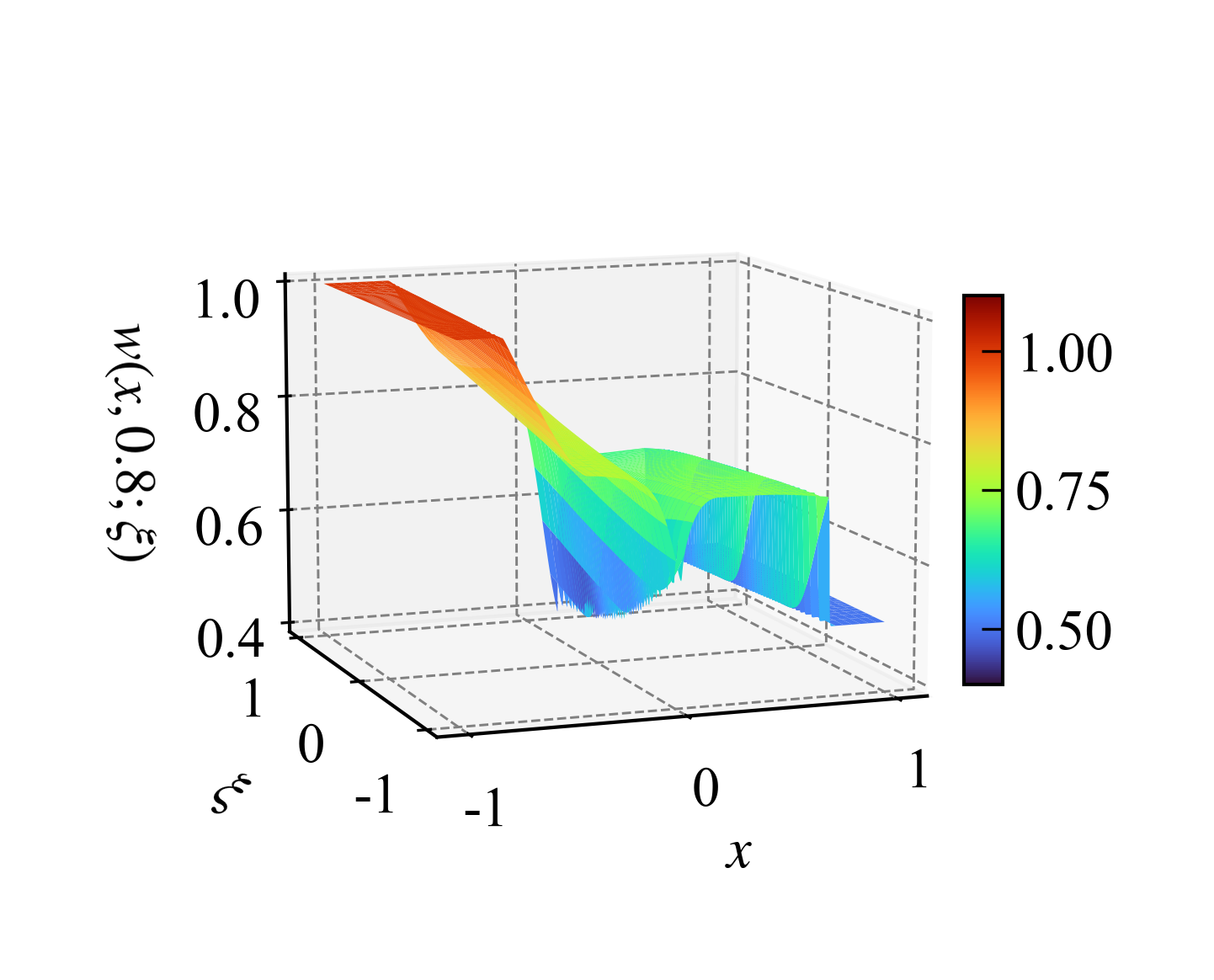}
\caption{}
\label{fig:3.4.1d}
\end{subfigure}
\hspace{0.05cm}
\centering
\begin{subfigure}[t]{0.26\textwidth}
\centering
\includegraphics[width=1.0\textwidth, trim=0.2cm 1.1cm 1.0cm 2.5cm, clip]{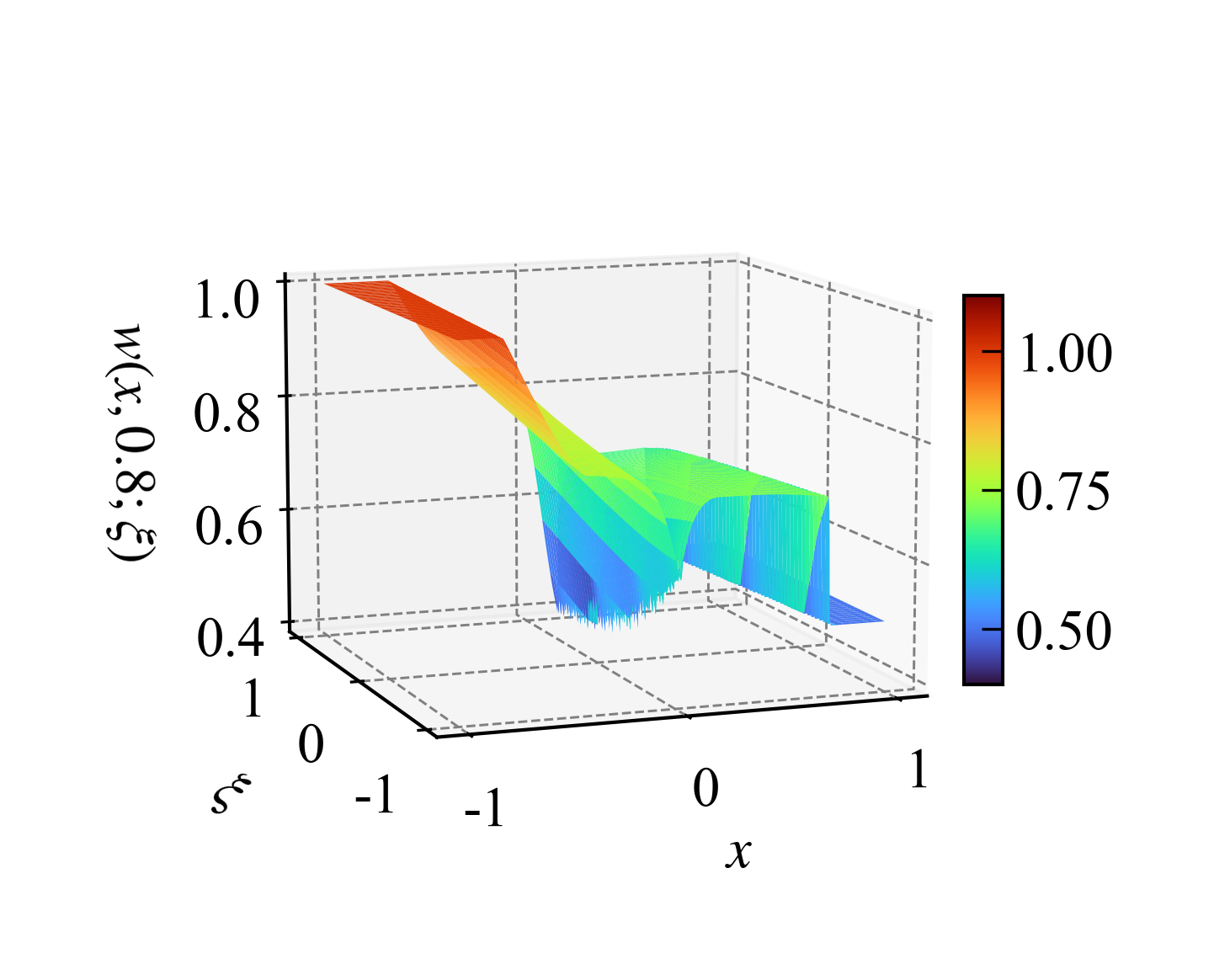}
\caption{}
\label{fig:3.4.1e}
\end{subfigure}
\hspace{0.05cm}
\centering
\begin{subfigure}[t]{0.26\textwidth}
\centering
\includegraphics[width=1.0\textwidth, trim=0.2cm 1.1cm 1.0cm 2.5cm, clip]{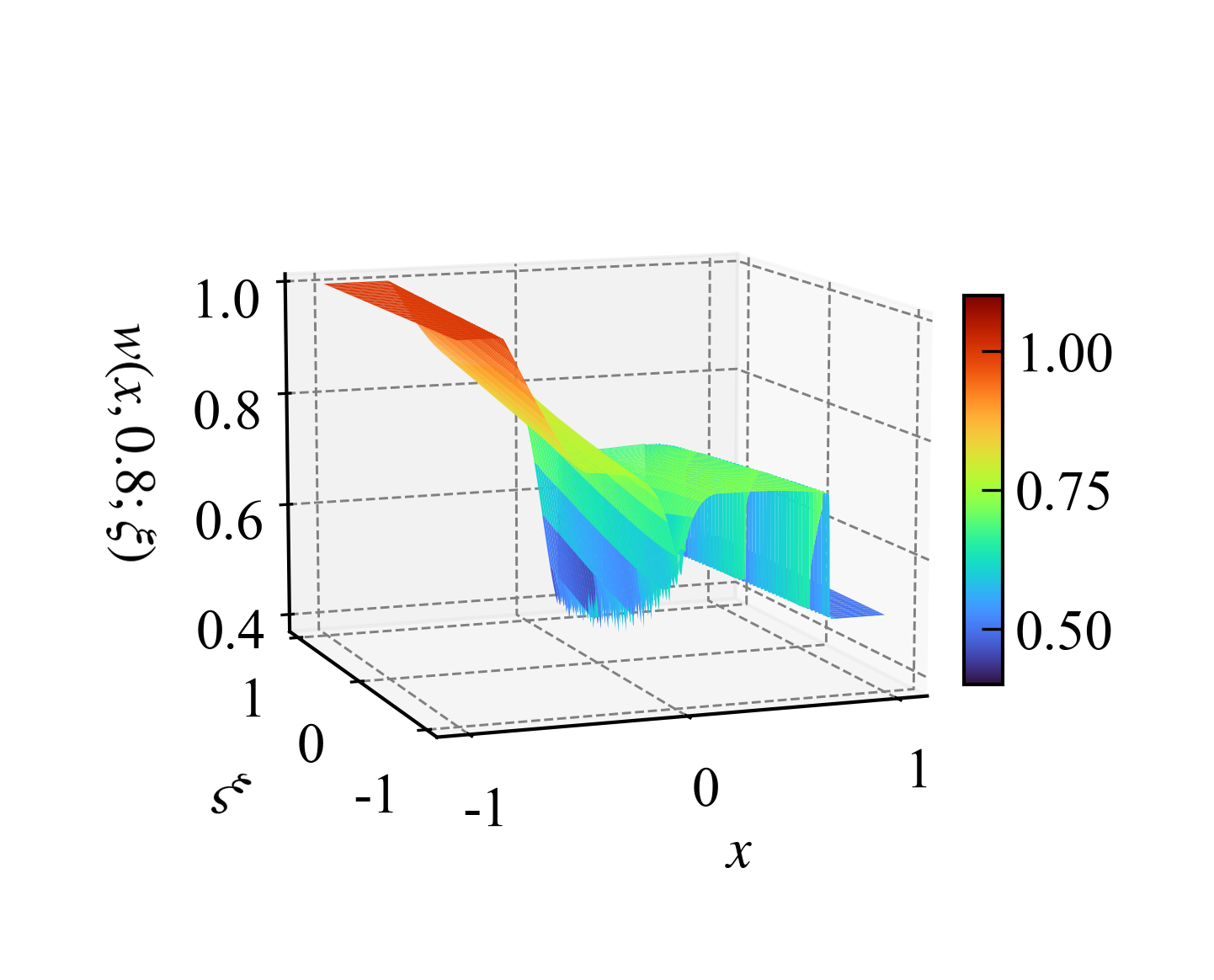}
\caption{}
\label{fig:3.4.1f}
\end{subfigure}
\captionsetup{justification=centering}
\caption{Example 4: $w(x,0.8;\xi)$ obtained in (\protect\subref{fig:3.4.1a})~simulations, (\protect\subref{fig:3.4.1b})~gPC, (\protect\subref{fig:3.4.1c})~interpolation B-spline, (\protect\subref{fig:3.4.1d})~approximation B-spline, (\protect\subref{fig:3.4.1e})~SP spline, and (\protect\subref{fig:3.4.1f})~CWENO.}
\label{fig:3.4.1}
\end{figure}

In \fref{fig:3.4.1}, we plot the water surface $w(x,0.8;\xi)$. \fref{fig:3.4.1a} shows the results obtained by solving the governing equations; the rest subfigures show the surface after application of the proposed methods. It is evident that the gPC solution and B-spline interpolation exhibit oscillations near the location of the discontinuity at $x\approx0.7$. When B-spline approximation, CWENO interpolation and SP splines are employed, these oscillations are suppressed. Note that the ``oscillations'' visible near $x\approx0$ appear in the solution.

\begin{figure}[ht!]
\centering
\begin{subfigure}[t]{0.26\textwidth}
\centering
\includegraphics[width=1.0\textwidth, trim=0.2cm 0.2cm 0.2cm 1.2cm, clip]{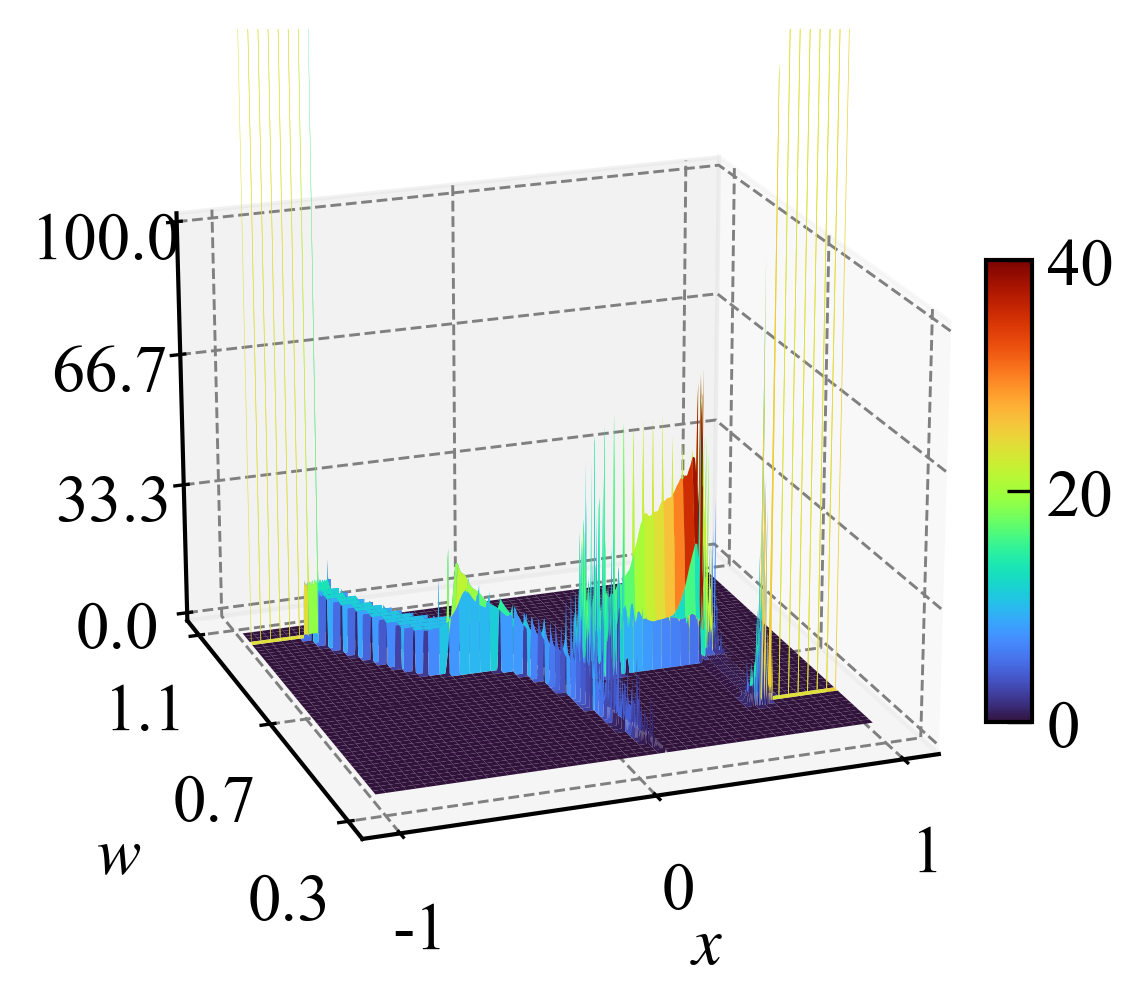}
\caption{}
\label{fig:3.4.2a}
\end{subfigure}
\centering
\begin{subfigure}[t]{0.26\textwidth}
\centering
\includegraphics[width=1.0\textwidth, trim=0.2cm 0.2cm 0.2cm 1.2cm, clip]{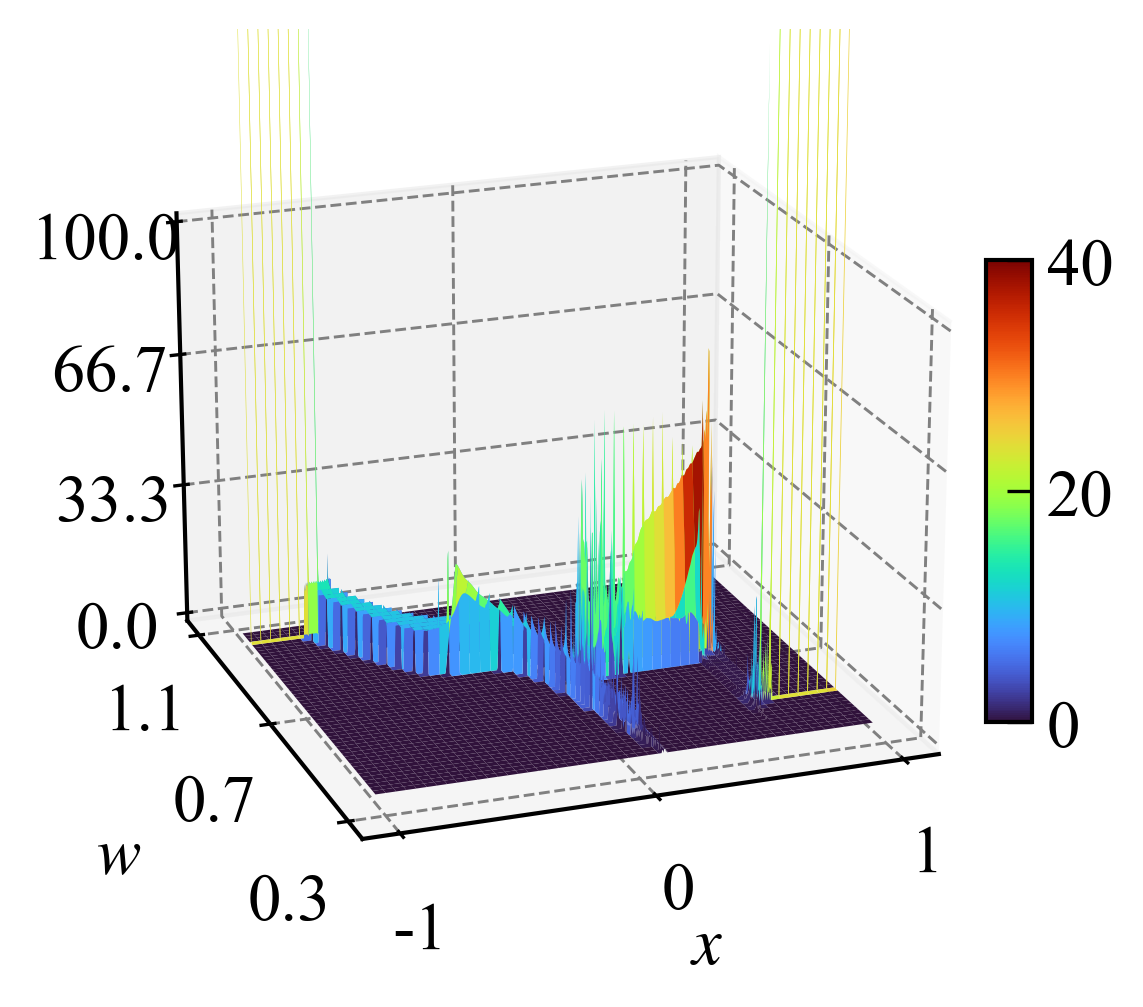}
\caption{}
\label{fig:3.4.2b}
\end{subfigure}
\centering
\begin{subfigure}[t]{0.26\textwidth}
\centering
\includegraphics[width=1.0\textwidth, trim=0.2cm 0.2cm 0.2cm 1.2cm, clip]{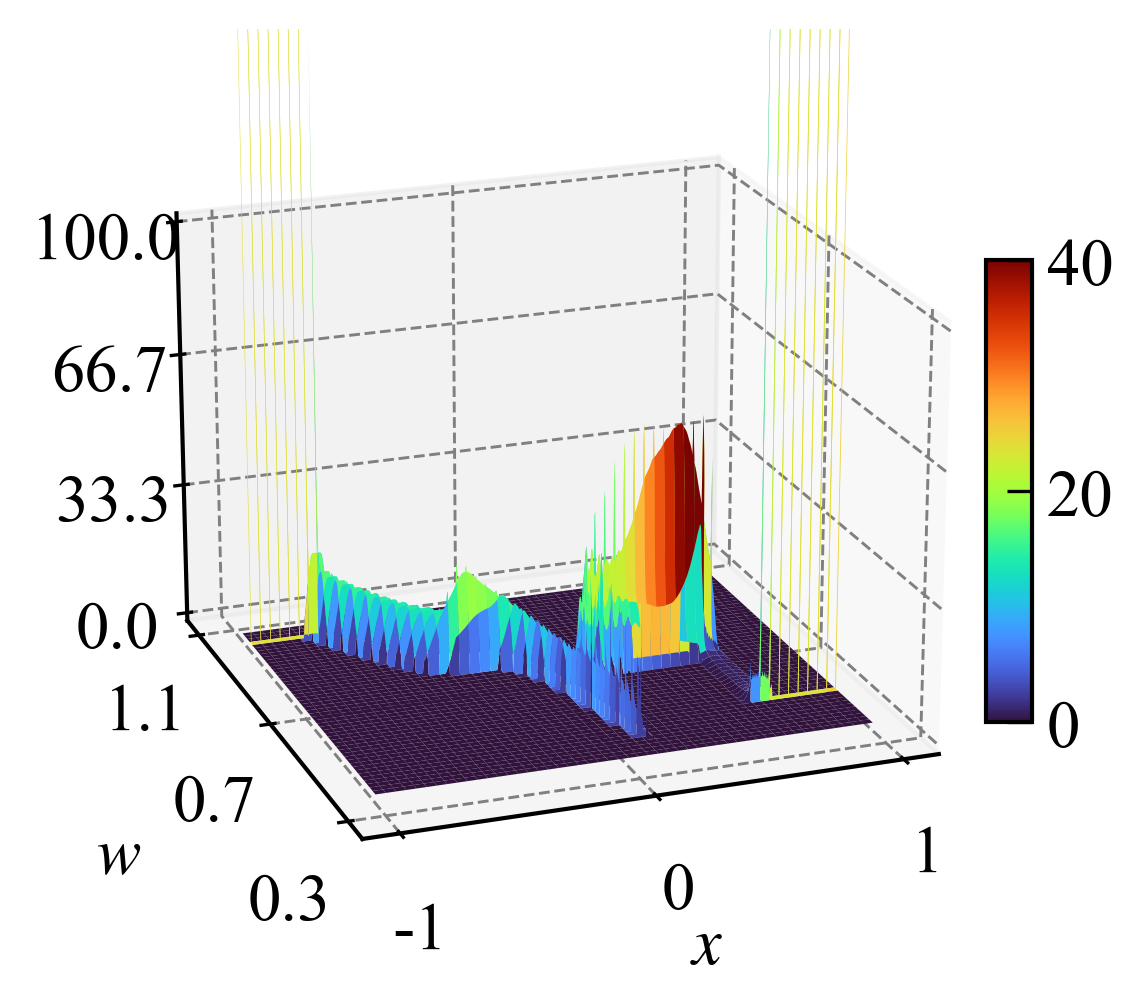}
\caption{}
\label{fig:3.4.2c}
\end{subfigure}
\centering
\begin{subfigure}[t]{0.26\textwidth}
\centering
\includegraphics[width=1.0\textwidth, trim=0.2cm 0.2cm 0.2cm 1.2cm, clip]{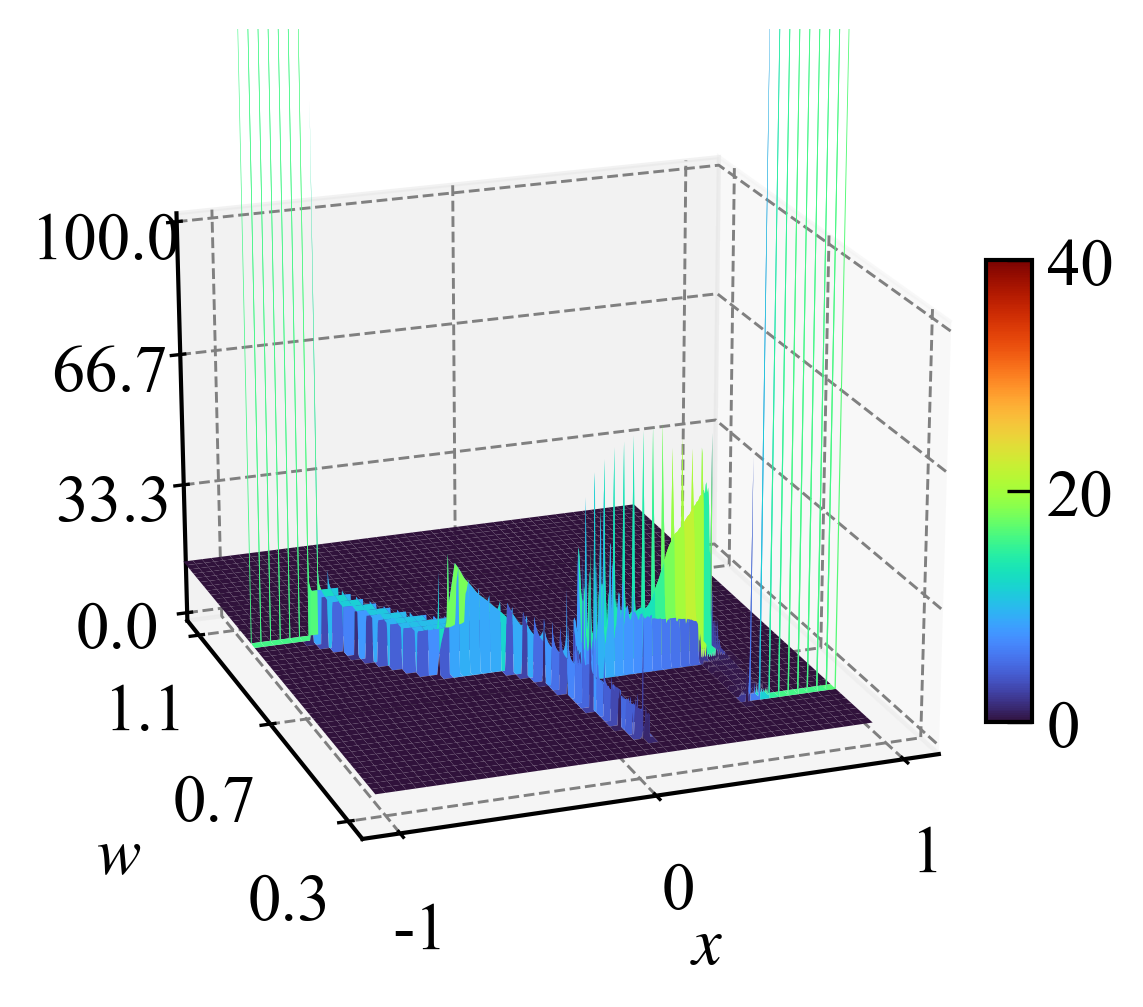}
\caption{}
\label{fig:3.4.2d}
\end{subfigure}
\centering
\begin{subfigure}[t]{0.26\textwidth}
\centering
\includegraphics[width=1.0\textwidth, trim=0.2cm 0.2cm 0.2cm 1.2cm, clip]{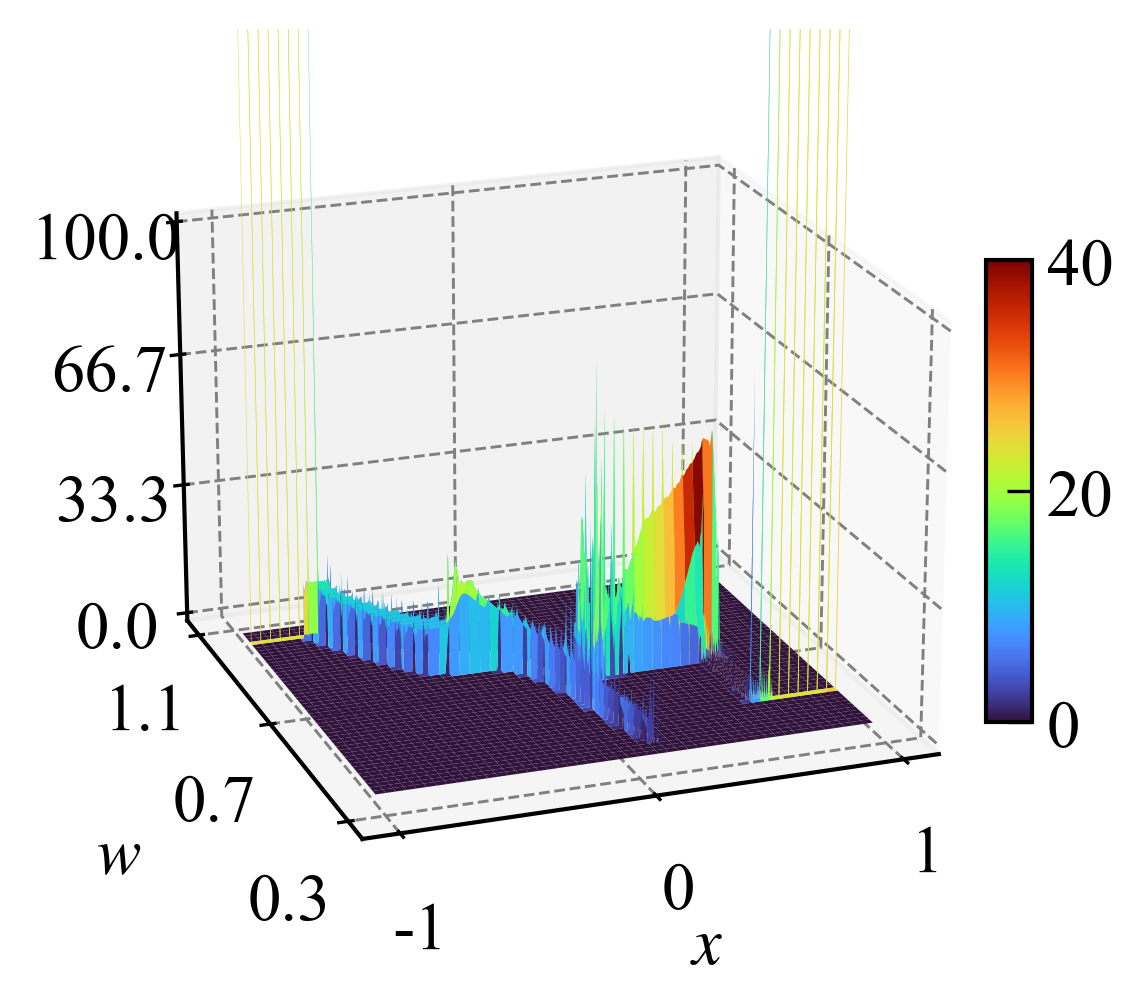}
\caption{}
\label{fig:3.4.2e}
\end{subfigure}
\centering
\begin{subfigure}[t]{0.30\textwidth}
\centering
\includegraphics[width=1.0\textwidth, trim=0.0cm 0.0cm 0.0cm 0.0cm, clip]{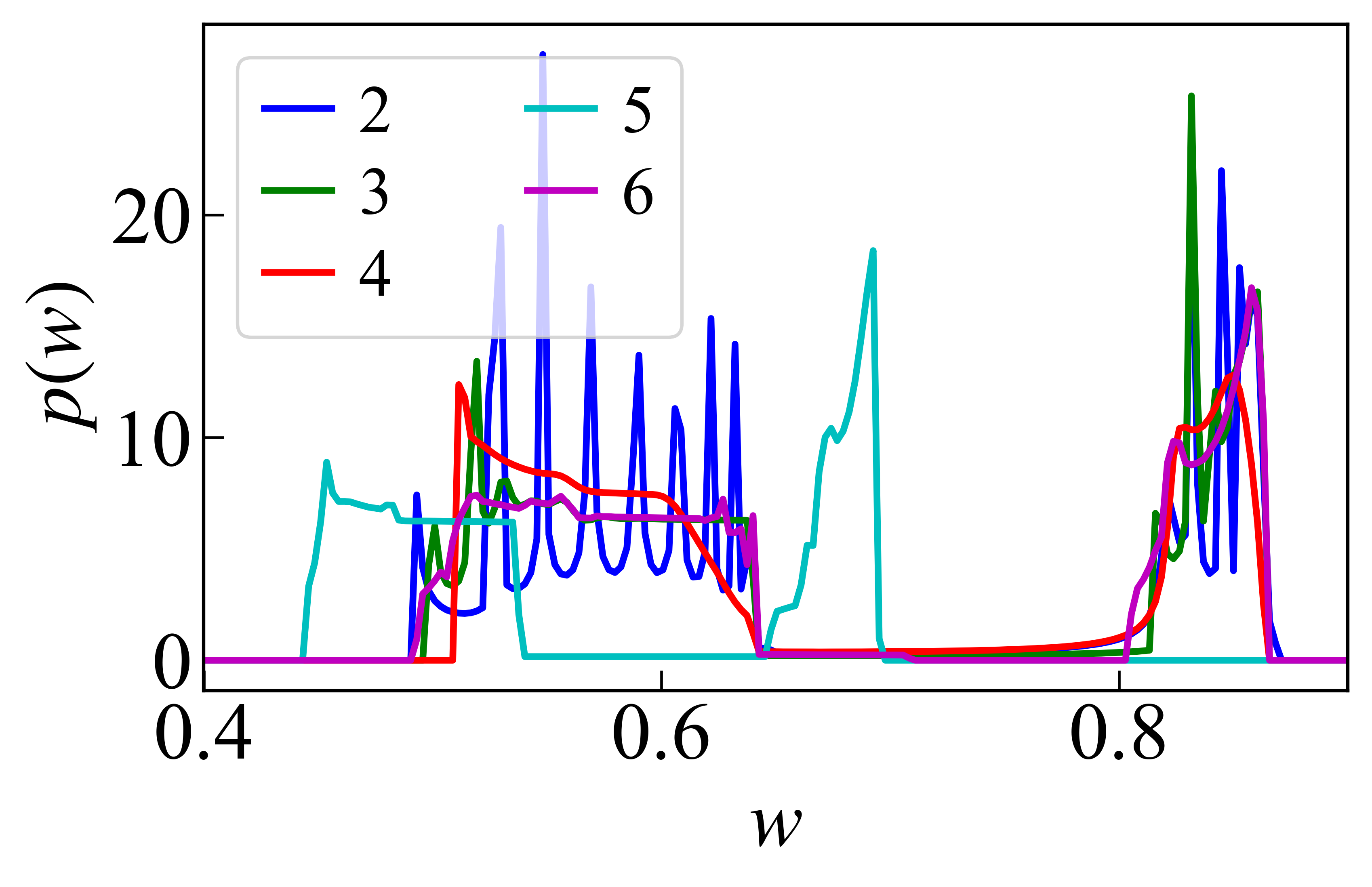}
\caption{}
\label{fig:3.4.2f}
\end{subfigure}
\captionsetup{justification=centering}
\caption{Example 4: PDF of $w(x,0.8;\xi)$ obtained using (\protect\subref{fig:3.4.2a})~gPC, (\protect\subref{fig:3.4.2b})~interpolation B-spline, (\protect\subref{fig:3.4.2c})~approximation B-spline, (\protect\subref{fig:3.4.2d})~SP spline, and (\protect\subref{fig:3.4.2e})~CWENO; (\protect\subref{fig:3.4.2f}) shows profile at $x=0.09$.}
\label{fig:3.4.2}
\end{figure}

The interpolated/approximated values are used to compute the PDF, shown in~\fref{fig:3.4.2}. All methods produce similar results, except for the approximation B-spline and SP spline. Additionally, gPC and interpolation B-spline exhibit more oscillations than CWENO interpolation, as illustrated in~\fref{fig:3.4.2f}, which shows a profile at $x=0.09$. 

\section{Conclusions}
This study presented a comparative analysis of several SC methods -- gPC, B-splines, SP splines, and CWENO interpolation. Each method’s performance was evaluated through a series of numerical examples. Our findings show that
\begin{enumerate}
\item For smooth solutions, gPC consistently provided accurate PDF reconstruction and exhibited the fastest convergence. However, it showed limitations in handling discontinuities, leading to oscillatory behavior and inaccurate PDF reconstruction.
\item Interpolation B-spline performed reliably for smooth data but, like gPC, encountered issues with discontinuous data due to oscillations. 
\item Approximation B-spline demonstrated limitations, particularly with its inherent smoothing properties, which impacted PDF accuracy in both smooth and discontinuous cases.
\item SP spline produced satisfactory results with smooth data but underperformed compared to CWENO interpolation when handling sharp gradients or discontinuities.
\item CWENO interpolation emerged as a robust option for handling both smooth and discontinuous data. While it required a sufficient number of points to maintain smoothness, it was the most effective at managing discontinuities, producing stable PDF approximations without the oscillations.
\end{enumerate}
In summary, CWENO interpolation shows significant potential as a versatile tool for SC, particularly in scenarios with discontinuities. Several problems in nuclear engineering involve discontinuous behavior due to abrupt transitions in physical states. Examples include water hammer events, steam explosions, and fuel-coolant interactions, among others. Future work might include advancing CWENO methods for SC to allow desired asymptotic behavior (e.g., non-negativity of water surface in shallow water equations) and exploring their applications in more complex multidimensional problems.



\bibliographystyle{mc2025}
\bibliography{main}


\end{document}